\documentclass{jhrs}
\usepackage{amsmath,epic,eepic,xy,url}
\xyoption{all}

\newtheorem{theorem}{Theorem}[section]
  \newtheorem{proposition}[theorem]{Proposition}
  \newtheorem{lemma}[theorem]{Lemma}
  \newtheorem{corollary}[theorem]{Corollary}

  \theoremstyle{definition}
  \newtheorem{definition}[theorem]{Definition}
  \newtheorem{example}[theorem]{Example}

  \theoremstyle{remark}
  \newtheorem{remark}[theorem]{Remark}
  
  \numberwithin{equation}{section}
\newcommand{\oil}{\mathopen{]}}    %% open interval left resp. right
\newcommand{\oir}{\mathclose{[}}   %% gives the good spacing.
\renewcommand{\epsilon}{\varepsilon}  %% nice epsilon and phi
\renewcommand{\phi}{\varphi}
\newcommand{\rest}[1]{\,|_{#1}}  %% restriction
\DeclareMathOperator{\id}{id}  %% identity mapping
\DeclareMathOperator{\Rep}{Rep_+}
\DeclareMathOperator{\Homeo}{Homeo_+}
\DeclareMathOperator{\Mon}{Mon_+}
\newcommand{\Pot}{\mathfrak{P}}  %% power set
\newcommand{\PotI}{\Pot_{[\,]}}  %% set of intervals
\DeclareMathOperator{\intr}{int}  %% interior (of a set)
\renewcommand*\le\leqslant
\renewcommand*\ge\geqslant

\begin{document}
\title{Reparametrizations of Continuous Paths}\thanks{The authors
  would like to thank the referee for many useful comments and in
  particular for pointing out an inaccuracy in the original proof of
  Theorem \ref{prop:tracehomeo}.}

\author{Ulrich Fahrenberg}
\email{uli@cs.aau.dk}
\homepage{www.cs.aau.dk/~uli}
\address{Department of Computer Science\\ Aalborg University\\
  Denmark\\ Fredrik Bajersvej 7B\\ DK-9220 Aalborg {\O}st}

\author{Martin Raussen}
\email{raussen@math.aau.dk}
\homepage{www.math.aau.dk/~raussen}
\address{Department of Mathematical Sciences\\ Aalborg University\\
    Denmark\\Fredrik Bajersvej 7G\\ DK-9220 Aalborg {\O}st}

\classification{55,68}

\keywords{path, regular path, reparametrization, reparametrization
  equivalence, trace, stop map, d-space, concurrency}

\begin{abstract}
  A reparametrization (of a continuous path) is given by a surjective
  weakly increasing self-map of the unit interval. We show that the
  monoid of reparametrizations (with respect to compositions) can be
  understood via ``stop-maps'' that allow to investigate compositions
  and factorizations, and we compare it to the distributive lattice of
  countable subsets of the unit interval. The results obtained are
  used to analyse the space of traces in a topological space, i.e.,
  the space of continuous paths up to reparametrization equivalence.
  This space is shown to be homeomorphic to the space of regular paths
  (without stops) up to increasing reparametrizations. Directed
  versions of the results are important in directed homotopy theory.
\end{abstract}

\received{June 9, 2006}
\revised{March 15, 2007}
\published{June 30, 2007}
\submitted{Ronnie Brown}

\volumeyear{2007}
\volumenumber{2}
\issuenumber{1}
\startpage{1}

\maketitle

\section{Introduction and Outline}

\subsection{Introduction}\label{ss:intro}

In elementary \emph{differential geometry}, the most basic objects studied
(after points perhaps) are \emph{paths}, i.e., \emph{differentiable} maps
$p:I\to \mathbf{R}^n$ defined on the closed interval $I=[0,1]$. Such a
path is called \emph{regular} if $p'(t)\neq \boldsymbol{0}$ for all $t\in
\oil 0,1\oir$. A \emph{reparametrization} of the unit interval $I$ is a
surjective differentiable map $\phi:I\to I$ with $\phi '(t)>0$ for all
$t\in \oil 0,1\oir$, i.e.~a (strictly increasing) self-diffeomorphism of
the unit interval.

Given a path $p:I\to \mathbf{R}^n$ and a reparametrization $\phi:I\to I$,
the paths $p$ and $p\circ \phi$ represent the same geometric object.  In
differential geometry one investigates equivalence classes (identifying
$p$ with $p\circ \phi$ for any reparametrization $\phi$) and their
invariants, like curvature and torsion.

Motivated by applications in \emph{concurrency theory}, a branch of
theoretical Computer Science trying to model and to understand the
coordination between many different processors working on a common task,
we are interested in \emph{continuous} paths $p:I\to X$ in more general
topological spaces up to more general reparametrizations $\phi:I\to I$.
When the \emph{state space} of a concurrent program is viewed as a
topological space (typically a cubical complex; cf.~\cite{FGR-ATC}),
``directed'' paths in that space respecting certain ``monotonicity''
properties correspond to \emph{executions}. A nice framework to handle
\emph{directed} topological spaces (with an eye to homotopy properties) is
the concept of a d-space proposed and investigated by Marco Grandis in
\cite{Grandis-DHTP}. Essentially, a topological space comes equipped with
a subset of preferred \emph{d-paths} in the set of all paths in $X$,
cf.~Definition \ref{def:d}. Note in particular, that the reverse of a
directed path in general is \emph{not} directed; the slogan is ``breaking
symmetries''.

We do not try to capture the \emph{quantitative} behaviour of executions,
corresponding to particular parametrizations of paths, but merely the
\emph{qualitative} behaviour, such as the order of shared resources used,
or the result of a computation. Hence the objects of study are paths up to
certain reparametrizations which
\begin{enumerate}
\item do not alter the \emph{image} of a path, and
\item do not alter the \emph{order of events}.
\end{enumerate}

We are thus interested in general paths in topological spaces, up to
\emph{surjective} reparametrizations $\phi: I\to I$ which are
\emph{increasing} (and thus continuous!---cf.~Lemma \ref{lem:char}), but
not necessarily strictly increasing. Two paths are considered to have
the same behaviour if they are \emph{reparametrization
  equivalent}, cf.~Definition \ref{def:repeq}.

To understand this equivalence relation, we have to investigate the
space of all reparametrizations which includes strange (e.g.~nowhere
differentiable) elements. Nevertheless, it enjoys remarkable
properties: It is a monoid, in which compositions and factorizations
can be completely analysed through an investigation of \emph{stop
  intervals} and of \emph{stop values}. The quotient space after
dividing out the self-homeomorphisms has nice algebraic lattice
properties.

A path is called \emph{regular} if it does not ``stop'';  and we are able
to show that the space of general paths \emph{modulo reparametrizations} is
homeomorphic to the space of \emph{regular} paths \emph{modulo increasing
  auto-homeomorphisms} of the interval. Hence to investigate properties of
the former, it suffices to consider the latter. This is one of the starting
points in the homotopy theoretical and categorical investigation of
invariants of d-spaces in \cite{Raussen:06}. Further possible areas of
application of the results (and of higher-dimensionsal generalisations
still to be investigated) include categorical homotopy theory
as in \cite{GJ:99}, categorified gauge theory as in \cite{BS:05} and
$n$-transport theory; cf.~the blog ``The $n$-category caf\'e'' at
\texttt{http://golem.ph.utexas.edu/category}.

This article does not build on any sophisticated machinery. Most of
the concepts and proofs can be understood with an undergraduate
mathematical background. There are certain parallels to the elementary
theory of distribution functions in probability theory,
cf.~e.g.~\cite{Loeve}. The flavour is nevertheless different, since
continuity (no jumps, i.e., surjectivity) is essential for us. For the
sake of completeness, we have chosen to include also elementary
results and their proofs (some of which may be well-known).

Marco Grandis has studied piecewise linear reparametrizations in
\cite{Grandis:05} for different purposes, but also in the framework of
``directed algebraic topology''.

\subsection{Basic definitions}

Let always $X$ denote a Hausdorff topological space and $I=[0,1]$ the unit
interval. The set of all (nondegenerate) closed subintervals of $I$ will
be denoted by $\PotI(I)=\{[a,b]\mid 0\le a<b\le 1\}$.  Let $p:I\to X$
denote a continuous map (a path), and remark that the pre-image
$p^{-1}(x)$ of any element $x\in X$ is a closed set.

\begin{definition}\label{def:stopreg}
  \begin{enumerate}
  \item An interval $J\in \PotI(I)$ is called a $p$-\emph{stop interval}
    if the restriction $p\rest J$ is constant and if $J$ is a maximal
    interval with that property.
  \item The set of all $p$-stop intervals will be denoted as $\Delta
    _p\subseteq \PotI(I)$. Remark that the intervals in $\Delta _p$ are
    disjoint and that $\Delta _p$ carries a natural total order.  We let
    $D_p:=\bigcup_{J\in \Delta_p} J\subset I$ denote the \emph{stop set}
    of $p$.
  \item A path $p:I\to X$ is called \emph{regular} if $\Delta
    _p=\emptyset$ or if $\Delta _p=\{I\}$ (no stop or
    constant).
  \item A continuous map $\phi:I\to I$ is called a
    \emph{reparametrization} if $\phi(0)=0, \phi (1)=1$ and if $\phi$ is
    \emph{increasing}, i.e.\ if $s\le t\in I$ implies $\phi (s)\le \phi
    (t)$.
  \end{enumerate}
\end{definition}

Remark that neither a regular path nor a reparametrization need be
injective.

\begin{definition}
  \label{def:repeq}
  Two paths $p, q: I\to X$ are called \emph{reparametrization equivalent}
  if there exist reparametrizations $\phi$, $\psi$ such that $p\circ \phi=
  q\circ \psi$.
\end{definition}

We will show later (Corollary~\ref{co:repeq-eq}) that
reparametrization equivalence is indeed an equivalence relation. As in
differential geometry, we are interested in equivalence classes of
paths modulo reparametrization equivalence. We call these equivalence
classes \emph{traces}\footnote{with a geometric meaning; the notion
  has nothing to do with algebraic traces.} in the space $X$. In
particular, we would like to know \emph{whether every trace can be
  represented by a regular path}. The (positive) answer to this
question in Proposition~\ref{prop:sur} is based on a closer look at
the space of reparametrizations of the unit interval.

\subsection{Outline of the article}

Section \ref{s:repar} contains a detailed study of reparametrizations
(in their own right) and characterizes their behaviour essentially by
an order-preserving bijection between the set of stop intervals and
the set of stop values (Definition \ref{def:stopreg} and Proposition
\ref{prop:repar-const}). This pattern analysis allows to study
compositions, and in particular, factorizations in the monoid of
reparametrizations from an algebraic point of view. In particular,
Proposition \ref{prop:Cbij} shows that the space of all
reparametrizations ``up to homeomorphisms'' is a distributive lattice
isomorphic to the lattice of countable subsets of the unit interval.

Section \ref{se:traces} investigates the \emph{space} of all paths in
a Hausdorff space up to reparametrization equivalence. The main result
(Theorem \ref{prop:tracehomeo}) states that two quotient spaces are in
fact homeomorphic: the orbit space arising from the action of the
group of all oriented homeomorphisms of the unit interval on the space
of \emph{regular} paths (with given end points, cf.~Definition
\ref{def:stopreg}) on the one side, and the space of \emph{all} paths
with given end points up to reparametrization equivalence
(Definition~\ref{def:repeq}); in particular, \emph{every trace can be
  represented by a regular path}. It might be a bit surprising that
the proof makes essential use of the results on factorizations of
reparametrizations from Section \ref{s:repar}.

The final Section \ref{s:dirtra} deals with spaces of \emph{directed}
traces (directed paths up to reparametrization equivalence) on a
d-space (cf.~Section \ref{ss:intro} and Definition \ref{def:d}).
Corollary \ref{cor:dsur} confirms that the result of Theorem
\ref{prop:tracehomeo} has an analogue for directed paths in
\emph{saturated} (cf.~Definition \ref{def:sat}) d-spaces. This result
is one of the starting points for the (categorical) investigations
into invariants of directed spaces in \cite{Raussen:06}. Furthermore,
it is shown how to relate reparametrization equivalence of directed
paths to \emph{thin dihomotopies}; this result is needed in the study
\cite{Fahrenberg:06} of directed squares (``two-dimensional paths'').

\section{Reparametrizations}\label{s:repar}

\subsection{Stop and move intervals, stop values, stop maps}

The following definitions (extending Definition~\ref{def:stopreg}) and
elementary results will mainly be used for reparametrizations. For the
sake of generality, we will state and prove them for general paths $p:
I\to X$ in a Hausdorff space $X$.

\begin{definition}\label{def:stop}
  \begin{enumerate}
  \item An element $c\in X$ is called a $p$-\emph{stop value} if there is
    a $p$-stop interval $J\in \Delta _p$ with $p(J)=\{c\}$. We let
    $C_p\subseteq X$ denote the set of all $p$-stop values.
  \item The map $p$ induces the $p$-\emph{stop map} $F_p:\Delta _p\to C_p$
    with $F_p(J)=c\Leftrightarrow p(J)=\{c\}.$
  \item An interval $J\in \PotI(I)$ is called a $p$-\emph{move
      interval} if it does not contain any $p$-stop interval and if it
    is maximal with that property.
  \item The set of all $p$-move intervals will be denoted $\Gamma
    _p\subseteq \PotI(I)$, a collection of disjoint closed intervals. We
    let $O_p:=\bigcup_{J\in \Gamma_p} \intr J\subseteq I$ denote the
    $p$-\emph{move set}.
  \end{enumerate}
\end{definition}

\begin{lemma}\label{lem:countable}
  For any path $p:I\to X$, the sets of $p$-stop intervals $\Delta _p$, of
  $p$-move intervals $\Gamma_P$ and of $p$-stop values $C_p$ are at most
  countable.
\end{lemma}

\begin{proof}
  The set $O_p$ and $\bigcup _{J\in \Delta_p}\intr J$ of
  \emph{interior points} in move, resp.~stop intervals are open
  subsets of $I$ and thus unions of at most \emph{countably} many
  maximal open intervals. Their closures constitute $\Gamma_p$,
  resp.~$\Delta _p$. The stop value set $C_p$ is at most countable as
  image of $\Delta _p$ under the $p$-stop map $F_p$.
\end{proof}

\begin{remark}
  This result is similar in spirit to the assertion (relevant for
  distribution functions in probability theory) that a nondecreasing
  function to an interval has at most countably many discontinuity points,
  cf.~e.g.~\cite[Sec.~11]{Loeve}.
\end{remark}

It is important to analyse the \emph{boundary} $\partial D_p$ of the
$p$-stop set: It can be decomposed as $\partial D_p=\partial _1D_p \cup
\partial _2D_p$ as follows:
\begin{itemize}
\item $\partial _1D_p=\partial D_p\cap D_p$ -- the set of all boundary
  points of intervals in $D_p$, an at most countable set;
\item $\partial _2D_p=\partial D_p\setminus D_p$ -- the set of all
  (honest) accumulation points of these boundary points. $\partial
  _2D_p$ can be uncountable; compare Ex.~\ref{ex:dense}.
\end{itemize}

The move set $O_p$ is the complement $O_p=I\setminus
\overline{D_p}\subset I$ of the closure of $D_p$. It
does occur that $O_p$ is empty; compare Ex.~\ref{ex:dense}.

The following elementary technical lemma concerning stop sets will
be needed in the proof of Proposition~\ref{prop:sur}.

\begin{lemma}\label{lem:invar}
  Let $p:I\to X$ denote a path and $U\subseteq X$ an open subspace.  Then
  $p^{-1}(U)$ is a union of \emph{(}at most\emph{)} countably many
  disjoint open intervals, and for any maximal open interval $\oil
  a,b\oir$ in $p^{-1}(U)$, $[a,c]\not\subseteq D_p$ and
  $[c,b]\not\subseteq D_p$ for every $c\in \oil a,b\oir$.
 \end{lemma}

\begin{proof}
  As an open subset of $I$, $p^{-1}(U)$ is a union of (at most)
  countably many disjoint open intervals. Let $\oil a,b\oir$ be one of
  these, and let $c\in \oil a,b\oir$. Then $p(c)\in U, p(a)\not\in U,
  p(b)\not\in U$. In particular, $p$ is neither constant on $[a,c]$
  nor on $[c,b]$.
\end{proof}

\subsection{Spaces of reparametrizations}\label{s:sparep}

Within the set of all self-maps of the unit interval $I$ fixing its
boundary points, we study the following subsets:

\begin{definition}
  \begin{itemize}
  \item $\Mon(I):=\{ \phi: I\to I\mid \phi\text{ increasing}, \phi(0)=0,
    \phi(1)=1 \}$;
  \item $\Rep(I):= \{ \phi\in \Mon(I)\mid \phi\text{ continuous }\}$ --
    \\the set of all increasing reparametrizations;
  \item $\Homeo(I) = \{ \rho\in \Rep(I)\mid \rho\text{ strictly increasing
    }\}$ -- \\the set of all increasing auto-homeomorphisms of the
    interval.
  \end{itemize}
\end{definition}

Note that $\Homeo(I)\subset \Rep(I)\subset \Mon(I)$. The compact-open
topology on the space of \emph{all} continuous maps $C(I,I)$ induces
topologies on the latter two spaces $\Rep(I)$ and $\Homeo(I)$.
Composition $\circ$ of maps turns $\Mon(I)$ into a monoid, $\Rep(I)$
into a topological monoid and $\Homeo(I)$ into a topological group
(consisting of the units in $\Rep(I)$).

All three mapping sets come equipped with a natural partial order: $\phi
\le \psi$ if and only if $\phi(t)\le \psi (t)$ for all $t\in I$, and they
form complete lattices with respect to $\le $. Least upper bounds, resp.\
greatest lower bounds are given by the $\max$, resp.\ $\min$ of the
functions involved:
\begin{equation*}
  (\phi \vee \psi)(t):=\max\{\phi(t),\psi (t)\}\qquad
  (\phi \wedge \psi)(t):=\min\{\phi(t),\psi(t)\}
\end{equation*}

\begin{lemma}\label{lem:convex}
  \begin{enumerate}
  \item All three sets $\Mon(I), \Rep(I), \Homeo(I)$ are convex. In
    particular, the latter two spaces are contractible.
  \item Any two reparametrizations $\phi, \psi\in \Rep(I)$ are
    d-homotopic (cf.~Definition \ref{def:d-hom} for the general
    definition), i.e.\ there exists a reparametrization $\phi, \psi\le
    \eta \in \Rep(I)$ and increasing paths $G,H: \vec{I}\to \Rep(I)$ with
    $G(0)=\phi, H(0)=\psi, G(1)= H(1)= \eta$.
\end{enumerate}
\end{lemma}

\begin{proof}
  \begin{enumerate}
  \item The sets are closed under convex combinations $(1-s)\phi
    +s\psi$.
  \item For $\eta= \phi \vee \psi$, define $G(s)=(1-s)\phi +s\eta$
    and $H(s)=(1-s)\psi +s\eta$.
  \end{enumerate}
\end{proof}

A characterization of the elements of $\Mon(I)$, $\Rep(I)$, and
$\Homeo(I)$ is achieved in the elementary
\begin{lemma}\label{lem:char}
  Let $\phi\in \Mon(I)$.
  \begin{enumerate}
  \item For every interval $J\subseteq I$, the pre-image
    $\phi^{-1}(J)\subseteq I$ is an interval, as well. In particular,
    $\phi^{-1}(a)$ is an interval (possibly degenerate) for every $a\in
    I$.
  \item $\phi\in \Rep(I)$ if and only if $\phi$ is \emph{surjective}.
  \item $\phi\in \Homeo(I)$ if and only if $\phi$ is \emph{bijective}.
  \end{enumerate}
\end{lemma}

\begin{proof}
  The only non-obvious statement is that surjectivity of $\phi \in
  \Mon(I)$ implies continuity; we show that the pre-image $\phi^{-1}(J)$
  of an open interval $J\subset I$ is open:

  Let $d\in \phi^{-1}(J)$ and $\phi (d)=c\in J$. Then there exist
  $\epsilon >0$ such that $[c-\epsilon, c+\epsilon]\subseteq J$
  and $d_1, d_2\in I$ such that $\phi(d_1)=c-\epsilon,
  \phi(d_2)=c+\epsilon$. Monotonicity implies:
  $\phi([d_1,d_2])\subseteq [c-\epsilon, c+\epsilon]$ and $d'\not\in
  [d_1,d_2]\Rightarrow \phi(d')\not\in \oil c-\epsilon,
  c+\epsilon\oir$.  Surjectivity implies:
  $\phi([d_1,d_2])=[c-\epsilon, c+\epsilon]\subseteq J$; hence $d$
  has an open neighbourhood in $\phi^{-1}(J)$.
\end{proof}

The following information about images of intervals under
reparametrizations is needed in the proof of Proposition \ref{prop:sur}:

\begin{lemma}\label{lem:pinterv}
  Let $a,b\in I$ and $\phi \in \Rep(I)$. Then
  \begin{enumerate}
  \item $\phi ([a,b])=[\phi (a),\phi (b)]$.
  \item $\oil \phi (a),\phi (b)\oir \subseteq \phi (\oil a,b\oir
    )\subseteq [\phi (a),\phi (b)]$.
  \item $\phi (\oil a,b\oir )\neq \oil \phi (a),\phi (b)\oir $ if and only
    if there is $c\in \oil a,b\oir$ such that $[a,c]\subseteq D_\phi$ or
    $[c,b]\subseteq D_\phi$.
  \end{enumerate}
\end{lemma}

\begin{proof}
  Only the last assertion requires proof. If $[a,c]\subseteq D_\phi$ for
  some $c\in \oil a, b\oir$, then $\phi(a)= \phi(c)\in \phi (\oil a,b\oir
  )$; similarly, if $[c,b]\subseteq D_\phi$, then $\phi(b)\in \phi (\oil
  a,b\oir )$. For the reverse direction, assume $\phi(a)\in \phi (\oil
  a,b\oir )$, and let $c\in \oil a, b\oir$ such that $\phi(a)= \phi(c)$.
  Then $\phi(a)\le \phi(t)\le \phi(c)= \phi(a)$ for any $t\in [a,c]$,
  hence $[a,c]\subseteq D_\phi$. The other implication is similar.
\end{proof}

The following result deals with the relative size of the
homeomorphisms within the reparametrizations. It will be needed in the
proof of the main result in Section~\ref{se:traces}.

\begin{lemma}\label{lem:dense}
  In the topology induced from the compact-open topology, both
  $\Homeo(I)$ and its complement are dense in $\Rep(I)$.
\end{lemma}

\begin{proof}
  The compact-open topology is induced by the supremum metric on the
  space $C(I,I)$ of all self-maps of the interval. Hence, for a given
  $\phi\in \Rep(I)$ and $n\in \mathbf{N}$, we need to construct
  $\rho\in \Homeo(I)$ such that $\| \phi-\rho\| \le \frac{1}{n}$:
  Choose $c_k, 0\le k\le n$, such that $c_0=0, c_n=1$ and
  $\phi(c_k)=\frac{k}{n}$; clearly $c_k$ is strictly increasing with
  $k$.  Hence the piecewise linear map $\rho$ given by $\rho
  (c_k)=\frac{k}{n}$ is contained in $\Homeo(I)$. Furthermore, for
  $x\in [c_k,c_{k+1}],k<n,$ we have $\frac{k}{n}\le \rho(x),\phi(x)\le
  \frac{k+1}{n}$, and thus $\| \phi-\rho\| \le \frac{1}{n}$.

  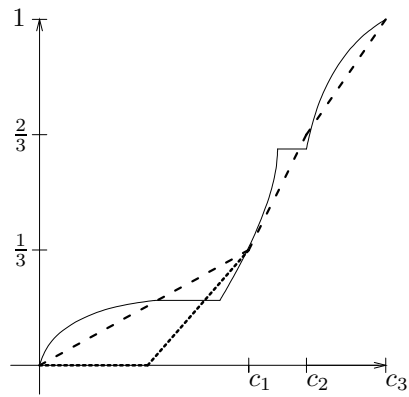
\begin{figure}
  \centering
  \setlength{\unitlength}{0.00050000in}
  \begingroup\makeatletter\ifx\SetFigFont\undefined%
  \gdef\SetFigFont#1#2#3#4#5{%
    \reset@font\fontsize{#1}{#2pt}%
    \fontfamily{#3}\fontseries{#4}\fontshape{#5}%
    \selectfont}%
  \fi\endgroup%
  {\renewcommand{\dashlinestretch}{30}
    \begin{picture}(3934,3949)(0,-10)
      \path(342.000,3792.000)(312.000,3912.000)(282.000,3792.000)
      \path(312,3912)(312,12)
      \path(12,312)(3837,312)
      \path(1512,987)(2187,987)
      \path(2787,2562)(3087,2562)
      \path(237,1512)(387,1512)\put(30,1442){$\frac{1}{3}$}
      \path(237,2712)(387,2712)\put(30,2642){$\frac{2}{3}$}
      \path(237,3912)(387,3912)\put(30,3842){$1$}
      \thicklines
      \dashline{90.000}(312,312)(2487,1512)(3087,2712)(3912,3912)
      \dottedline{68}(312,312)(1437,312)(2487,1512)
      \thinlines
      \path(2487,312)(2487,237)(2487,387)\put(2487,100){$c_1$}
      \path(3087,387)(3087,237)\put(3087,100){$c_2$}
      \path(3087,312)(3912,312)
      \path(3792.000,282.000)(3912.000,312.000)(3792.000,342.000)
      \path(3912,387)(3912,237)\put(3912,100){$c_3$}
      \path(312,312)(312,313)(313,316)
      (314,321)(317,328)(320,338)
      (324,350)(330,366)(336,383)
      (344,403)(353,423)(364,446)
      (375,469)(389,492)(404,517)
      (422,541)(442,567)(465,593)
      (491,620)(521,647)(556,675)
      (595,704)(639,733)(687,762)
      (731,785)(776,807)(820,827)
      (865,845)(908,862)(951,876)
      (992,890)(1033,901)(1073,912)
      (1113,922)(1152,930)(1190,938)
      (1227,945)(1264,952)(1299,958)
      (1333,963)(1365,968)(1394,972)
      (1421,976)(1444,979)(1464,981)
      (1480,983)(1493,985)(1502,986)
      (1507,986)(1511,987)(1512,987)
      \path(2187,987)(2188,988)(2189,991)
      (2192,995)(2196,1002)(2202,1013)
      (2210,1026)(2219,1042)(2231,1062)
      (2244,1085)(2258,1110)(2274,1138)
      (2291,1168)(2309,1200)(2328,1233)
      (2347,1268)(2366,1304)(2386,1340)
      (2406,1378)(2426,1417)(2446,1457)
      (2467,1499)(2488,1542)(2509,1586)
      (2530,1632)(2552,1680)(2573,1730)
      (2595,1782)(2616,1834)(2637,1887)
      (2659,1947)(2679,2003)(2697,2057)
      (2712,2106)(2725,2151)(2736,2193)
      (2746,2231)(2754,2267)(2760,2300)
      (2766,2332)(2771,2361)(2775,2388)
      (2778,2414)(2780,2439)(2782,2461)
      (2784,2482)(2785,2500)(2786,2516)
      (2786,2530)(2787,2541)(2787,2549)
      (2787,2555)(2787,2559)(2787,2561)(2787,2562)
      \path(3087,2562)(3087,2563)(3088,2566)
      (3089,2572)(3091,2580)(3093,2592)
      (3096,2607)(3100,2625)(3105,2647)
      (3111,2673)(3117,2700)(3124,2731)
      (3132,2763)(3141,2796)(3150,2831)
      (3161,2867)(3172,2905)(3184,2943)
      (3197,2981)(3211,3021)(3226,3063)
      (3244,3105)(3262,3149)(3283,3194)
      (3306,3241)(3331,3289)(3358,3338)
      (3387,3387)(3419,3438)(3452,3486)
      (3484,3530)(3515,3570)(3545,3607)
      (3575,3640)(3603,3671)(3631,3699)
      (3658,3724)(3684,3748)(3710,3770)
      (3735,3790)(3760,3809)(3783,3826)
      (3805,3842)(3826,3857)(3845,3869)
      (3862,3881)(3877,3890)(3888,3897)
      (3898,3903)(3904,3907)(3909,3910)
      (3911,3911)(3912,3912)
    \end{picture}}
  \caption{ \label{fig:dense}Reparametrization $\phi$ (full line),
    homeomorphism $\rho$ (broken line), and reparametrization $\psi$
    (dotted line).}
\end{figure}

For the same $\phi\in \Rep(I)$ and the same definition for $c_0$ and $c_1$
as above, let $\psi\in \Rep(I)\setminus \Homeo(I)$ be given by
$\psi(x)=\phi(x), x\ge c_1$; on the interval $[0,c_1]$, we let $\psi$ be
the piecewise linear map with $\psi(0)=\psi(\frac{c_1}{2})=0$ and
$\psi(c_1)=\phi(c_1)$. See also Figure~\ref{fig:dense}.
\end{proof}

\subsection{Classification of reparametrizations}

In the following, we are mainly interested in an investigation of the
algebraic monoid structure on $\Rep(I)$ induced by \emph{composition}
$\circ$ of maps. Note that there is another structure on the sets (spaces)
$\Mon(I)$, $\Rep(I)$, and $\Homeo(I)$, induced by \emph{concatenation} of
paths
\begin{equation}\label{eq:concat}
  (\phi, \psi)\mapsto \phi* \psi;\quad (\phi* \psi)(t)=
  \begin{cases}
    \phi(2t) &\text{for $t\le \frac12$} \\
    \psi(2t-1) &\text{for $t> \frac12$}
  \end{cases}
\end{equation}
This composition does not induce a monoidal structure on these sets, as
concatenation is not associative and does not have units ``on the nose''.

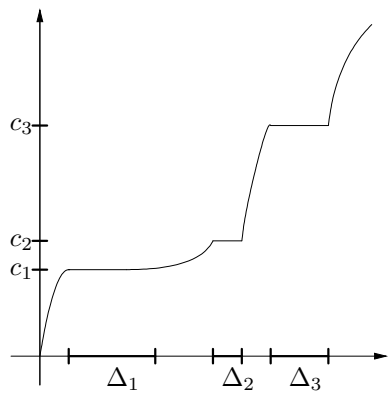
\begin{figure}
  \centering
  \setlength{\unitlength}{0.0005in}
  \begingroup\makeatletter\ifx\SetFigFont\undefined%
  \gdef\SetFigFont#1#2#3#4#5{%
    \reset@font\fontsize{#1}{#2pt}%
    \fontfamily{#3}\fontseries{#4}\fontshape{#5}%
    \selectfont}%
  \fi\endgroup%
  \renewcommand{\dashlinestretch}{30}
  \begin{picture}(3924,3939)(0,-10)
    \blacken\path(342.000,3792.000)(312.000,3912.000)(282.000,3792.000)(342.000,3792.000)
    \path(312,3912)(312,12) \path(12,312)(3912,312)
    \blacken\path(3792.000,282.000)(3912.000,312.000)(3792.000,342.000)(3792.000,282.000)
    \drawline(612,1212)(612,1212) \path(612,1212)(1212,1212)
    \path(2112,1512)(2412,1512) \path(2712,2712)(3312,2712)
    \path(312,12)(312,3912) \thicklines \path(612,312)(1512,312)
    \path(2112,312)(2412,312) \path(2712,312)(3312,312)
    \path(237,1212)(387,1212) \path(237,1512)(387,1512)
    \path(237,2712)(387,2712) \path(612,387)(612,237)
    \path(1512,387)(1512,237) \path(2112,387)(2112,237)
    \path(2412,387)(2412,237) \path(2712,387)(2712,237)
    \path(3312,387)(3312,237) \thinlines \path(312,312)(312,313)(313,317)
    (314,323)(315,332)(318,345) (320,362)(324,382)(328,405)
    (333,432)(338,460)(344,491) (350,523)(356,556)(362,590)
    (369,624)(376,659)(384,694) (392,730)(400,766)(409,802)
    (419,839)(429,877)(439,914) (451,951)(462,987)(477,1029)
    (491,1066)(504,1097)(517,1123) (528,1143)(538,1160)(548,1173)
    (558,1183)(566,1192)(575,1198) (582,1203)(589,1206)(595,1209)
    (601,1210)(605,1211)(608,1212) (610,1212)(611,1212)(612,1212)
    \path(1212,1212)(1214,1212)(1217,1212)
    (1223,1212)(1233,1212)(1246,1212) (1263,1213)(1283,1213)(1307,1214)
    (1333,1214)(1361,1215)(1390,1217) (1421,1218)(1452,1221)(1485,1223)
    (1518,1226)(1552,1230)(1586,1235) (1622,1241)(1659,1247)(1696,1255)
    (1735,1264)(1774,1275)(1812,1287) (1860,1304)(1902,1322)(1938,1340)
    (1969,1358)(1994,1375)(2016,1392) (2035,1409)(2051,1426)(2065,1442)
    (2077,1457)(2087,1471)(2095,1484) (2102,1494)(2107,1502)(2110,1507)
    (2111,1511)(2112,1512) \path(2412,1512)(2412,1513)(2412,1516)
    (2413,1522)(2414,1530)(2415,1541) (2417,1555)(2419,1571)(2421,1590)
    (2424,1610)(2427,1633)(2430,1658) (2434,1684)(2439,1713)(2444,1745)
    (2451,1780)(2458,1819)(2467,1863) (2476,1911)(2487,1962)(2497,2006)
    (2506,2049)(2515,2088)(2523,2125) (2531,2157)(2537,2184)(2543,2208)
    (2548,2228)(2552,2246)(2556,2261) (2559,2274)(2562,2287)(2565,2300)
    (2568,2313)(2572,2327)(2576,2343) (2581,2361)(2587,2382)(2593,2407)
    (2601,2434)(2609,2464)(2618,2496) (2627,2530)(2637,2562)(2644,2585)
    (2650,2605)(2656,2623)(2662,2640) (2667,2654)(2672,2666)(2676,2676)
    (2680,2685)(2683,2693)(2686,2699) (2689,2704)(2692,2708)(2694,2711)
    (2696,2714)(2698,2716)(2700,2717) (2702,2717)(2704,2718)(2705,2718)
    (2706,2717)(2707,2717)(2709,2716) (2709,2715)(2710,2715)(2711,2714)
    (2711,2713)(2712,2712) \path(3312,2712)(3312,2713)(3313,2717)
    (3313,2722)(3314,2731)(3316,2743) (3318,2759)(3320,2777)(3323,2799)
    (3327,2822)(3331,2848)(3335,2876) (3341,2905)(3346,2935)(3353,2967)
    (3360,2999)(3368,3032)(3377,3067) (3387,3104)(3399,3142)(3412,3182)
    (3427,3225)(3444,3268)(3462,3312) (3483,3359)(3504,3403)(3525,3443)
    (3545,3480)(3565,3513)(3584,3544) (3603,3572)(3621,3597)(3639,3621)
    (3656,3643)(3673,3664)(3689,3683) (3704,3700)(3718,3715)(3730,3729)
    (3740,3740)(3749,3748)(3755,3755) (3759,3759)(3761,3761)(3762,3762)

    \put(0,1162){$c_1$} \put(0,1462){$c_2$} \put(0,2662){$c_3$}
    \put(1000,0){$\Delta_1$} \put(2200,0){$\Delta_2$}
    \put(2900,0){$\Delta_3$}
  \end{picture}
  \caption{Stop intervals and stop values}
  \label{fig:cdelta}
\end{figure}

We wish to describe a reparametrization $\phi \in \Rep(I)$ by its
$\phi$-stop map $F_{\phi}:\Delta _{\phi}\to C_{\phi}$ illustrated in
Figure~\ref{fig:cdelta} and by the restriction of $\phi$ to its
$\phi$-move set $O_{\phi}\subseteq I$; cf.~Definition \ref{def:stop}.

\subsubsection{All countable sets in the interval are stop value sets}

Lemma \ref{lem:countable} tells us that the $\phi$-stop value set
$C_\phi\subset I$ of a reparametrization $\phi$ is an at most countable
ordered subset of $I$. Remark also that automorphisms $\rho \in \Homeo
(I)$ are characterized by the properties $\Delta_\rho=C_\rho=\emptyset$,
resp.\ $O_\rho= I$.

Which (countable) subsets of the unit interval can be realized as
$\phi$-stop sets of some $\phi\in\Rep(I)$? It is easy to construct
(piecewise linear) reparametrizations with a \emph{finite} set of stop
values. Rather surprisingly, this construction can be extended to
arbitrary (at most) \emph{countable} sets of stop values:

\begin{lemma}
  \label{lem:count}
  For every countable set $C\subset I$, there is a reparametrization
  $\phi\in \Rep(I)$ with $C_{\phi}=C$.
\end{lemma}

\begin{proof}
  Let $C=\{c_1,c_2,\dots \}\subset I$ denote an injective enumeration of the
  countable set $C$. We shall first construct a uniformly convergent
  sequence of piecewise linear maps $\phi _n\in \Rep(I), n\ge 0$ with
  $C(\phi _n)=\{c_1,\dots ,c_n\}$ and thus $\Delta _{\phi _n}=\{ \phi
  _n^{-1}(c_i)|1\le i\le n\} $. Let $[x_i^-,x_i^+]=\phi _n^{-1}(c_i),
  1\le i$; moreover, $x_0^-=0, x_0^+=1$.

  We start with $\phi_0= \id_I$.  Inductively, assume $\phi _n$ given
  as above. Among the $x_j^{\pm}, 1\le j\le n,$ choose $x_i^+,x_k^-$
  such that $c_{n+1}\in \phi( \oil x_i^+,x_k^-\oir )$ and such that
  the restriction of $\phi_n$ on that interval is strictly increasing
  (and linear). The map $\phi _{n+1}$ will differ from $\phi _n$ only
  on (the interior of) that subinterval $[c_i,c_k]$. The linear map on
  that interval is replaced by a piecewise linear map, which comes in
  three pieces.  The middle one takes the constant value $c_{n+1}$ on
  a subinterval $[x_{n+1}^-,x_{n+1}^+]$. On the left and right
  subinterval, we connect linearly to the values $c_i,c_k$ on the
  boundaries. See also Figure~\ref{fig:insert}.

  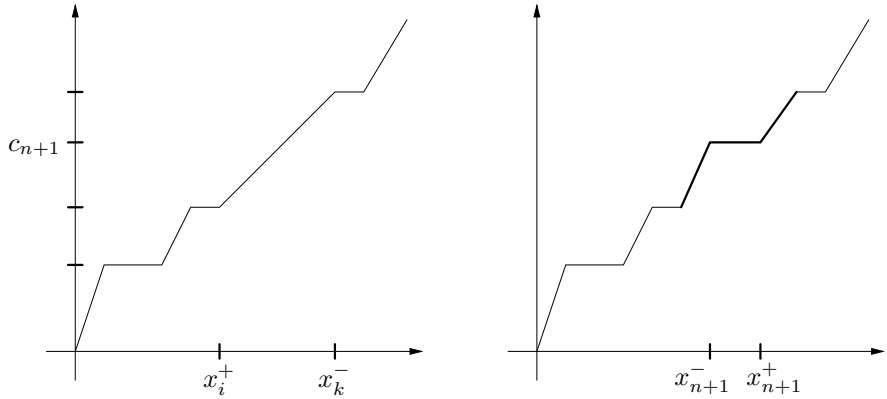
\begin{figure}
    \centering
\setlength{\unitlength}{0.00050000in}
\begingroup\makeatletter\ifx\SetFigFont\undefined%
\gdef\SetFigFont#1#2#3#4#5{%
  \reset@font\fontsize{#1}{#2pt}%
  \fontfamily{#3}\fontseries{#4}\fontshape{#5}%
  \selectfont}%
\fi\endgroup%
{\renewcommand{\dashlinestretch}{30}
\begin{picture}(9669,4144)(0,-10)
\blacken\path(1287.000,3997.000)(1257.000,4117.000)(1227.000,3997.000)(1287.000,3997.000)
\path(1257,4117)(1257,217)
\path(957,517)(4857,517)
\blacken\path(4737.000,487.000)(4857.000,517.000)(4737.000,547.000)(4737.000,487.000)
\blacken\path(6087.000,3997.000)(6057.000,4117.000)(6027.000,3997.000)(6087.000,3997.000)
\path(6057,4117)(6057,217)
\path(5757,517)(9657,517)
\blacken\path(9537.000,487.000)(9657.000,517.000)(9537.000,547.000)(9537.000,487.000)
\path(1257,517)(1557,1417)(2157,1417)
	(2457,2017)(2757,2017)(3957,3217)
	(4257,3217)(4707,3967)
\thicklines
\path(1182,1417)(1332,1417)
\path(1182,2017)(1332,2017)
\path(1182,3217)(1332,3217)
\thinlines
\path(8757,3217)(9057,3217)(9507,3967)
\path(6057,517)(6357,1417)(6957,1417)
	(7257,2017)(7557,2017)
\thicklines
\path(7557,2017)(7857,2692)(8382,2692)(8757,3217)
\path(1182,2692)(1332,2692)
\path(2757,592)(2757,442)
\path(3957,592)(3957,442)
\path(7857,592)(7857,442)
\path(8382,592)(8382,442)
\put(1107,2542){\makebox(0,0)[rb]{{\SetFigFont{10}{12.0}{\rmdefault}{\mddefault}{\updefault}$c_{n+1}$}}}
\put(2757,67){\makebox(0,0)[b]{{\SetFigFont{10}{12.0}{\familydefault}{\mddefault}{\updefault}$x_i^+$}}}
\put(3957,67){\makebox(0,0)[b]{{\SetFigFont{10}{12.0}{\familydefault}{\mddefault}{\updefault}$x_k^-$}}}
\put(7782,67){\makebox(0,0)[b]{{\SetFigFont{10}{12.0}{\familydefault}{\mddefault}{\updefault}$x_{n+1}^-$}}}
\put(8532,67){\makebox(0,0)[b]{{\SetFigFont{10}{12.0}{\familydefault}{\mddefault}{\updefault}$x_{n+1}^+$}}}
\end{picture}
}
\caption{ \label{fig:insert}Inserting the stop value $c_{n+1}$}
  \end{figure}

  The interval $[x_{n+1}^-,x_{n+1}^+]$ is chosen so small that $||\phi
  _{n+1}-\phi _n||_{\infty}<\frac{1}{2^n}$ ensuring uniform convergence of
  the maps $\phi _n$ to a continuous map $\phi \in \Rep(I)$. For this map
  $\phi$, we have $\Delta_{\phi}= \{ [x_n^-,x_n^+]|n>0\}$ and $C_{\phi}=C$.
\end{proof}

\begin{example}\label{ex:dense}
  If one chooses a \emph{dense} countable subset $C\subset I$, e.g.,
  $C=\mathbf{Q}\cap I$, then $\phi$ cannot be injective on \emph{any}
  non-trivial interval; hence $O_{\phi}=\emptyset$ and
  $\overline{D_{\phi}}=I$.  The (uncountable) complement $I\setminus
  D_{\phi}=\partial _2D_{\phi}$ does not contain any non-trivial
  interval.
\end{example}

\subsubsection{Classification}

What are the essential data to describe a reparametrization in terms
of stop maps and move sets?

\begin{proposition}\label{prop:desc}
  Let $\phi\in \Rep(I)$ denote a reparametrization.
  \begin{enumerate}
  \item The reparametrization $\phi$ induces an \emph{order-preserving
      bijection} $F_{\phi}:\Delta _{\phi}\to C_{\phi}$. The restriction
    $\phi\rest{J}: J\to \phi(J)$ to every move interval $J\in \Gamma_\phi$
    (Definition~\ref{def:stop}) is an (increasing) homeomorphism onto its
    image.
  \item The restriction $\phi\rest{\overline{D_\phi}}:
    \overline{D_\phi}\to \overline{C_\phi}$ of $\phi$ to
    $\overline{D_{\phi}}$ is \emph{onto}.
  \item Two reparametrizations $\phi, \psi\in \Rep(I)$ with $\Delta
    _\phi=\Delta_\psi, C_\phi=C_\psi, F_\phi=F_\psi$ agree on
    $\overline{D_{\phi}}$; if, moreover, $\phi\rest {O_\phi}=\psi\rest
    {O_\psi}$, then $\phi$ and $\psi$ agree on all of $I$.
  \end{enumerate}
\end{proposition}

\begin{proof}
  \begin{enumerate}
  \item The first statement is obvious from the definitions. For the
    second, note that $J\cap D_{\phi}=\partial J$ (or empty) for every such
    interval $J$.
  \item Every element $b\in \overline{C_{\phi}}$ is the limit of a
    monotone sequence of elements in $C_{\phi}$ which is the image of
    a monotone and bounded sequence of elements in $D_{\phi}$; the
    limit of such a sequence exists and maps to $b$ under
    $\phi$.
  \item By definition, $\phi$ and $\psi$ agree on $D_{\phi}$; by
    continuity, they have to agree on its closure $\overline{D_\phi}$,
    as well. The last statement is obvious.
  \end{enumerate}
\end{proof}

A reparametrization $\phi\in \Rep(I)$ is thus \emph{uniquely}
characterized by its stop map $F_{\phi}:\Delta _{\phi}\to C_{\phi}$ and by
a (fitting) collection of homeomorphisms $\phi|_{\overline{J}}:
\overline{J}\to \phi(\overline{J}),\; J\in \Gamma_{\phi}$. Now we ask
which conditions an ``abstract'' stop map has to satisfy in order to arise
from a genuine reparametrization. We start with the following data:

\begin{itemize}
\item $\Delta\subseteq \PotI(I)$ denotes an (at most) countable subset of
  \emph{disjoint closed} intervals -- with a natural total order.
\item $C\subseteq I$ denotes a subset with the same cardinality as
  $\Delta$.
\item $F:\Delta\to C$ denotes an order-preserving bijection.
\end{itemize}

Let $\Delta_-, \Delta_+ \subseteq I$ denote the set of lower, resp.~upper
boundaries of intervals in $\Delta$. Define $D:=\bigcup _{J\in
  \Delta}\Delta\subset I$.

Let $O=I\setminus \overline{D}$. Since $O$ is open, it is a
disjoint union $O=\bigcup _{J\in \Gamma}J$ of maximal open intervals
indexed by an (at most) countable set $\Gamma$ -- possibly empty.

For every map $G:\Delta\to C$ we define a map $\phi_G:D\to C$ by
$\phi_G(t)=G(J)\Leftrightarrow t\in J$. If $F$ is order-preserving,
then $\phi_F$ is increasing. Moreover:

\begin{proposition}\label{prop:repar-const}
  \begin{enumerate}
  \item A reparametrization $\phi \in \Rep(I)$ satisfies the following
    for every pair of strictly monotonely converging sequences
    $x_n\!\uparrow\! x$, $y_n\!\downarrow\! y$ for which $x, x_n\in
    (\Delta_\phi)_+$ and $y, y_n\in (\Delta_\phi)_-$:
    \begin{align}
      x=y & \quad\Rightarrow\quad \lim \phi(x_n)=\lim \phi(y_n),\label{eqa:f.1}\\
      x<y & \quad\Rightarrow\quad \lim \phi(x_n)<\lim \phi(y_n),\label{eqa:f.2}\\
      x=1 & \quad\Rightarrow\quad \lim \phi(x_n)=1,\label{eqa:f.3}\\
      y=0 & \quad\Rightarrow\quad \lim \phi(y_n)=0,\label{eqa:f.4}
    \end{align}
  \item For every order preserving bijection $F:\Delta\to C$ with
    $\phi_F$ satisfying $(1)$--$(4)$ above for every pair of strictly
    monotonely converging sequences $x_n\!\uparrow\! x$,
    $y_n\!\downarrow\! y$ for which $x, x_n\in (\Delta_\phi)_+$ and
    $y, y_n\in (\Delta_\phi)_-$, there exists a reparametrization
    $\psi \in \Rep (I)$ with $\Delta_{\psi}= \Delta, C_{\psi}=C$ and
    $F_{\psi}=F$. The set of all such reparametrizations is in
    one-to-one correspondence with $\prod _{\Gamma} \Homeo (I)$.
  \end{enumerate}
\end{proposition}

\begin{proof}
  \begin{enumerate}
  \item By continuity, $\lim \phi (x_n)=\phi (x)$ and $\lim \phi
    (y_n)=\phi (y)$; this settles all but (\ref{eqa:f.2}). Suppose
    $\phi (x)=\phi (y)$ in (2). Then $[x,y]$ is contained in a stop
    interval, hence $x\not\in \Delta_+$ and $y\not\in \Delta_-$.
  \item First, we extend $\phi_F$ to $\overline{D}$: there is a
    \emph{unique continuous} (and increasing!) extension of $\phi_F$
    from $D$ to $\overline{D}$: $\lim_{x_n\uparrow x} \phi_F(x_n)$
    exists and is independent of the sequence $x_n$ by monotonicity
    and agrees with $\lim_{y_n\downarrow x} \phi_F(y_n)$ by
    condition~\eqref{eqa:f.1} of Proposition~\ref{prop:repar-const}.
    Moreover, we let $\phi_F(0)=0$ and $\phi_F(1)=1$, in accordance
    with (3) and (4) above.

    Let $J=\oil a_-^J,a_+^J\oir\in \Gamma$ denote a maximal open
    interval.  Its boundary points $a_-^J, a_+^J$ are contained in
    $\partial D$ unless possibly if $a_-^J=0$ and/or $a_+^J=1$, in
    which case we are covered by \eqref{eqa:f.3} and/or
    \eqref{eqa:f.4} above.  In conclusion, $\phi_F$ is defined on
    $\partial J$. Moreover, $\phi_F(a_-^J)<\phi_F(a_+^J)$, since $F$
    is order preserving and injective and because of condition
    \eqref{eqa:f.2} above.

    Hence, every collection of strictly increasing homeomorphisms
    between $[a_-^J,a_+^J]$ and $[\phi_F(a_-^J),\phi_F(a_+^J)]$ --
    preserving endpoints -- extends $\phi_F$ to a continuous
    increasing map $\psi:I\to I$ with $\Delta_{\psi}= \Delta,
    C_{\psi}=C$ and $F_{\psi}=F$. The set of all collections of such
    homeomorphisms is easily seen to be in one-to-one correspondence
    with $\prod _{\Gamma} \Homeo (I)$.
  \end{enumerate}
\end{proof}

\subsection{Compositions and Factorizations}
\label{se:cf}

We shall now investigate the behaviour of $\Rep(I)$ under composition
and factorization in view of the description and classification from
Proposition~\ref{prop:repar-const} above. We need to introduce the
following notation: For a (continuous) map $\psi: I\to I$, let $\psi
_*,(\psi^{-1})^*:\PotI(I)\to \PotI(I)$ denote the maps induced on
subintervals: $\psi _*(J)=\psi (J)\in \PotI(I),
(\psi^{-1})^*(J)=\psi^{-1}(J)$.

\subsubsection{Composition of reparametrizations}

The results below follow easily from the definitions of
stop-intervals, stop-values and stop-maps:

\begin{lemma}\label{lem:CD}
  Let $\phi, \psi\in \Rep(I)$ denote reparametrizations with associated
  stop maps $F_{\phi}:\Delta _{\phi}\to C_{\phi}, F_{\psi}:\Delta
  _{\psi}\to C_{\psi}$. Then
  \begin{enumerate}
  \item $\Delta_{\phi \circ \psi}=\{J\in\Delta_{\psi}\mid
    F_{\psi}(J)\not\in D_{\phi}\}\cup (\psi^{-1})^*(\Delta_{\phi}) $,
  \item $C_{\phi \circ \psi}=\phi (C_{\psi})\cup C_{\phi}$,
  \item $F_{\phi \circ \psi}: \Delta _{\phi \circ \psi}\to C_{\phi \circ
      \psi}$ is given by $F_{\phi \circ \psi}(J)=
    \begin{cases}
      \phi (F_{\psi}(J)), & J\in \Delta_{\psi}\\
      F_{\phi}(\psi _*(J)), & J\in (\psi^{-1})^*(\Delta_{\phi}).
    \end{cases}
    $
  \end{enumerate}
  \qed
\end{lemma}

\begin{corollary}\label{cor:CD}
  Let $\phi,\psi \in \Rep(I)$ as in Lemma \ref{lem:CD}.  If $\psi \in
  \Homeo(I)$, resp.~$\phi\in \Homeo(I)$, then
  \begin{enumerate}
  \item $\Delta_{\phi \circ \psi}=(\psi^{-1})^*(\Delta_{\phi}) $,
    resp.~$\Delta_{\phi \circ \psi}=\Delta _\psi$,
  \item $C_{\phi \circ \psi}=C_{\phi}$, resp.~$C_{\phi \circ
      \psi}=\phi(C_\psi)$,
  \item $F_{\phi \circ \psi}=F_\phi\circ \psi _*:\psi
    _*^{-1}(\Delta_{\phi})\to C_{\phi}$, resp.~$F_{\phi \circ
      \psi}=\phi\circ F_\psi: \Delta _\psi\to \phi(C_\psi)$.
  \end{enumerate}
  \qed
\end{corollary}

\subsubsection{Factorizations of reparametrizations}

Also factorizations can be studied effectively using stop-data of the
reparametrizations involved. It turns out that the following result on
factorizations on the \emph{right} will be an essential tool in
Section~\ref{se:traces}:
\begin{proposition}
  \label{prop:lift}
  Let $\eta ,\phi\in \Rep(I)$ denote reparametrizations.
  \begin{enumerate}
  \item There exists a lift $\psi \in \Rep(I)$ in the diagram
    \begin{equation}
      \label{eq:lift-cp}
      \xymatrix{%
        & I\ar[d]^\phi \\
        I\ar[r]^{\eta}\ar@{-->}[ur]^{\psi} & I
      }
    \end{equation}
    if and only if $C_\phi\subseteq C_{\eta}$.
  \item If $C_\phi\subseteq C_{\eta}$ and $C$ is any (at most)
    countable set with $\phi^{-1}(C_{\eta}\setminus C_\phi)\subseteq
    C\subseteq \phi^{-1}(C_{\eta}\setminus C_\phi)\cup D_\phi$, then
    there exists such a lift $\psi \in \Rep(I)$ with $C_{\psi}=C$.\\
    In particular, if $C_\phi=C_{\eta}$, there exist a lift $\psi
    \in \Homeo(I)$.
  \item Assume $C_\phi\subseteq C_{\eta}$ and let $\Delta
    _1:=F_\eta ^{-1}(C_\phi)\subseteq \Delta_\eta$. Then the space
    of all lifts $\{\psi \mid \eta =\phi\circ \psi\}\subseteq
    \Rep(I)$ is in one-to-one-correspondence with $\prod _{L\in
      \Delta_1}\Rep(L)=\prod _{\Delta_1}\Rep(I)$.
  \end{enumerate}
\end{proposition}

\begin{proof}
  The ``only if'' part of 1.~ follows immediately from Lemma
  \ref{lem:CD}.2.  For the ``if'' part, we analyse first the set-theoretic
  requirements to a lift $\psi $ on relevant subintervals. To this end,
  decompose $\Delta_\eta=:\Delta_1\sqcup \Delta_2$ with $\Delta
  _1:=F_\eta ^{-1}(C_\phi)$ and $\Delta
  _2=F_\eta^{-1}(C_\eta\setminus C_\phi)$, and $D_\eta =D_1\sqcup
  D_2$ with $D_1:=\bigcup_{J\in \Delta _1}$ and $D_2=\bigcup_{J\in \Delta
    _2}$.  We construct a lift $\psi: I=D_1\cup D_2\cup O_\eta\to I$ by
  considering each of these three subsets of $I$:

  Remark that $\Delta _2$ necessarily has to be a subset of $\Delta _\psi$
  and that for $J\in \Delta _2$, $F_\psi (J)$ has to be the unique element
  of $\phi^{-1}(F_\eta (J))$.

  On any move interval $K\in \Gamma_\eta$ (cf.~Definition
  \ref{def:stop}.4), the restriction $\eta \rest{K}: K\to \eta (K)$ is
  an increasing homeomorphism; in particular, $\eta (K)\cap C_\eta$
  consists at most of the two boundary points. Hence, the restriction
  $\phi \rest {\phi ^{-1}\eta (K)}: \phi ^{-1}\eta (K)\to \eta (K)$
  is also an increasing homeomorphism, since $\eta (K)\cap
  C_\phi\subseteq \eta (K)\cap C_\eta$ again consists at most of the
  two boundary points. The restriction of $\psi$ to $K$ has to be defined
  as $\psi\rest K=(\phi\rest{\phi^{-1}\eta K} )^{-1}\circ
  \eta\rest{K}$; it is onto $\phi^{-1}\eta K$.

  On any interval $L\in \Delta _1$, the restriction of $\psi $ to $L$ can
  be defined as \emph{any} increasing continuous map $\psi\rest L: L\to
  F_{\phi}^{-1}(F_\eta(L))\in \Delta_\phi$ respecting the boundary
  points.

  The map $\psi :\vec{I}\to \vec{I} $ thus defined altogether is by definition a
  lift, it is increasing and surjective. Lemma \ref{lem:char}.2 settles 1.

  The only freedom in the construction of a lift $\psi$ is the choice of
  increasing continuous maps on the intervals $L\in \Delta _1$ (with given
  end points). As in Lemma~\ref{lem:count}, we can construct the set of
  stop values (on $D_1$) to be any countable subset of $D_\phi$. To these
  one has of course to add the pre-images of the stop values in
  $C_{\eta}\setminus C_\phi$. This settles 2. and 3.
\end{proof}

We will also make use of the following result on factorizations on the
\emph{left}.

\begin{proposition}\label{prop:lift-left}
  Let $\eta ,\phi\in \Rep(I)$ denote reparametrizations.
  \begin{enumerate}
  \item There exists a factorization with $\psi \in \Rep(I)$ in the
    diagram
    \begin{equation}
      \label{eq:lift-left}
      \xymatrix{%
        I\ar[d]_\phi\ar[r]^{\eta} & I\\
        I\ar@{-->}[ur]_{\psi} &
      }
    \end{equation}
    if and only if there exists a map $i_{\phi\eta}:\Delta
    _\phi\to\Delta_\eta$ such that $J\subseteq i_{\phi\eta}(J)$
    for every $J\in \Delta _\phi$. $(\Delta_\phi$ is a
    \emph{refinement} of $\Delta_\eta)$.
  \item If it exists, the factor $\psi\in \Rep(I)$ is uniquely
    determined and satisfies
    \begin{itemize}
    \item $C_\psi=C_\eta\setminus \{F_\eta (J)\mid J\in \Delta_\eta
      \cap \Delta_\phi\}$,
    \item $\Delta _\psi=\{ \phi(K)\mid K\in \Delta_\eta \setminus
      \Delta_\phi\}$,
    \item if $K\in \Delta_\psi$, then $F_\psi (K)$ is the unique element
      of $\eta (\phi^{-1}(K)),\, K\in \Delta_\psi$.
    \end{itemize}
  \end{enumerate}
\end{proposition}

\begin{proof}
  A lift $\psi$ as in (\ref{eq:lift-left}) has to satisfy $\psi
  (x):=\eta(\phi^{-1}(x))$. It is well-defined (and then unique and
  increasing) if and only if the condition of
  Proposition~\ref{prop:lift-left} is satisfied. Since $\eta $ is onto,
  $\psi$ is onto as well, and thus continuous by Lemma \ref{lem:char}.2.
  The description of the invariants of $\psi$ follows by inspection.
\end{proof}

\subsection{The algebra of reparametrizations up to homeomorphisms}

Consider the group action $\Rep(I)\times \Homeo(I)\to \Rep(I)$ given
by composition on the right. An element in the quotient space
$\Rep(I)/_{\Homeo(I)}$ preserves the set of stop values, whereas the
exact distribution of stop intervals over the interval is factored
out. Using the factorization tools from Section~\ref{se:cf} above, this
intuition will be made more formal in Propositions \ref{prop:Cbij} and
\ref{prop:distlat} below.

Consider the preorder on $\Rep(I)$ (different from the one considered in
Section~\ref{s:sparep}) given by $\phi \le \psi \Leftrightarrow \exists\;
\eta\in \Rep(I): \psi = \phi\circ \eta$ ($\Leftrightarrow C_\phi\subseteq
C_\psi$ by Proposition \ref{prop:lift}.1).  This preorder factors to yield
a partial order on the quotient $\Rep(I)/_{\Homeo(I)}$ since $\psi =\phi
\circ \eta = (\phi \circ \rho)\circ (\rho ^{-1}\circ\eta)$ for $\rho \in
\Homeo (I)$.

Moreover, let us consider the set $\Pot_c(I)$ of \emph{countable} subsets
of $I$ with the partial order given by inclusion.

\begin{proposition}\label{prop:Cbij}
  The map $C:\Rep(I)/_{\Homeo(I)}\to \Pot_c(I)$ given by
  $C(\phi)=C_\phi$ is an order-preserving bijection.
\end{proposition}

\begin{proof}
  By Corollary \ref{cor:CD}.2, the map $C$ is well-defined; by Lemma
  \ref{lem:CD}, it is order-preserving, and by Lemma \ref{lem:count},
  it is surjective. Given two reparametrizations with the same set of
  stop values, Proposition~\ref{prop:lift} shows that one can
  construct a lift $\psi$ from one into the other that is a
  homeomorphism ($C_\psi =\emptyset $); as a consequence, $C$ is
  also injective.
\end{proof}

\begin{proposition}
  \label{prop:surjeq}
  For every $\phi_1, \phi_2\in \Rep(I)$, there exist $\psi_1, \psi_2\in
  \Rep(I)$ completing the diagram
  \begin{equation*}
    \xymatrix{%
      I \ar@{-->}[r]^{\psi_1} \ar@{-->}[d]_{\psi_2} & I \ar[d]^{\phi_1} \\ I
      \ar[r]_{\phi_2} & I
    }
  \end{equation*}
  with $C_{\phi_1\circ\psi_1}=C_{\phi_2\circ\psi_2}=C_{\phi_1}\cup
  C_{\phi_2}.$
\end{proposition}

\begin{proof}
  Using Lemma \ref{lem:count}, construct $\psi _1\in \Rep(I)$ with
  $C_{\psi _1}=\phi_1^{-1}(C_{\phi_2}\setminus C_{\phi _1})$ and hence
  $C_{\phi_1\circ \psi_1}=C_{\phi _1}\cup C_{\phi _2}$ (cf.~Lemma
  \ref{lem:CD}.2). Using Proposition~\ref{prop:lift}, construct a lift
  $\psi _2$ in the diagram
  \begin{equation*}
    \xymatrix{%
      & I\ar[d]^{\phi_2} \\
      I\ar[r]_{\phi_1\circ \psi_1}\ar@{-->}[ur]^{\psi_2} & I }
  \end{equation*}
  with $C_{\psi_2}=\phi_2^{-1}(C_{\phi_1}\setminus C_{\phi _2})$, i.e.,
  without introducing superfluous extra stop values.
\end{proof}

\begin{remark}
  In general, it is not possible to complete the dual diagram
  \begin{equation*}
    \xymatrix{%
      I \ar[r]^{\phi_1} \ar[d]_{\phi_2} & I \ar@{-->}[d]^{\psi_1} \\ I
      \ar@{-->}[r]_{\psi_2} & I }
  \end{equation*}
  since $D_{\psi_1\circ \phi_1}=D_{\psi_2\circ \phi_2}\supseteq
  D_{\phi_1}\cup D_{\phi_2}$. The latter set might be the entire
  interval $I$ which is impossible for a reparametrization.
\end{remark}

Is there a natural way to construct from two reparametrizations a third
one (a common factor) with a set of stop values that is just the
\emph{intersection} of the sets of stop values of the given ones? In order
to have the stop intervals appear in proper order, it is necessary to
modify one of the reparametrizations by a homeomorphism first (which is
not a problem if one works in the quotient $\Rep(I)/_{\Homeo(I)}$!)
\begin{proposition}\label{prop:gcd}
  For every $\phi_1, \phi_2\in \Rep(I)$, there exist $\rho\in\Homeo(I)$,
  $\psi_1, \psi_2, \phi \in \Rep(I)$ completing the diagram
  \begin{equation*}
    \xymatrix{I\ar[r]^{\phi
        _1}\ar@{-->}[d]_{\rho}\ar@{-->}[ddr]^{\phi} & I\\
      I\ar[d]_{\phi_2}& \\ I & I\ar@{-->}[l]^{\psi_2}\ar@{-->}[uu]_{\psi_1} }
  \end{equation*}
  and such that $C_\phi=C_{\phi_1}\cap C_{\phi_2}$.
\end{proposition}

\begin{proof}
  Define $\Delta_{\phi}:=F_{\phi_1}^{-1}(C_{\phi_1}\cap
  C_{\phi_2})\subseteq \Delta_{\phi_1}, C_\phi:=C_{\phi_1}\cap C_{\phi_2} $
  and define $F_{\phi}: \Delta_{\phi}\to C_{\phi}$ as the restriction
  of $F_{\phi_1}$.  By Proposition~\ref{prop:repar-const}, there
  exists $\phi\in \Rep(I)$ with $F_{\phi}$ as its stop map. By
  Proposition~\ref{prop:lift-left}, there exists a lift $\psi_1\in
  \Rep(I)$ in the right triangle of the diagram above.

  In general, the reparametrization $\phi$ constructed above does not
  factor over $\phi_2$ immediately,
  cf.~Proposition \ref{prop:lift-left}. We need a ''correction''
  homeomorphism $\rho\in \Homeo(I)$ whose restriction to $D_{\phi}$
  fits into
  \begin{equation*}
    \xymatrix{D_{\phi}\ar@{-->}[r]^{\rho}\ar[d]_{F_{\phi}} & D_{\phi
        _2}\ar[d]^{F_{\phi _2}}\\ C_{\phi}\ar[r]^{\subseteq}& C_{\phi_2}.}
  \end{equation*}
  On intervals $J\in \Delta_\phi$, $\rho$ can be chosen as the
  increasing linear map sending $J$ onto $F_{\phi_2}^{-1}(F_\phi (J))$
  and then extended from $D_\phi$ to $I$ as a homeomorphism as in the
  proof of Proposition~\ref{prop:repar-const}.  The condition of Proposition
  \ref{prop:lift-left} is now satisfied to guarantee a lift of $\phi$
  over $\phi_2\circ \rho$ in the left triangle of the diagram above.
\end{proof}

Using the bijection from Proposition~\ref{prop:Cbij}, one may
introduce binary operations on the quotient $\Rep(I)/_{\Homeo(I)}$ in
a purely algebraic manner, i.e., one may pull back the operations
given by set union and intersection on $\Pot_c(I)$. The results of
this section allow us to give these operations an intrinsic meaning in
terms of reparametrizations. Using the notation from
Proposition~\ref{prop:surjeq}, an operation $\vee $ (``least common
multiple'') is defined by $[\phi_1]\vee [\phi _2]:=[\phi _1\circ
\psi_1]=[\phi _2\circ \psi_2]$. Likewise, $[\phi_1]\wedge [\phi _2]$
can be represented by the reparametrization $\phi$ from
Proposition~\ref{prop:gcd}. Altogether we obtain:

\begin{proposition}\label{prop:distlat}
  The operations
  \begin{equation*}
    \vee, \wedge: \Rep(I)/_{\Homeo(I)}\times \Rep(I)/_{\Homeo(I)}
    \to \Rep(I)/_{\Homeo(I)}
  \end{equation*}
  turn $\Rep(I)/_{\Homeo(I)}$ into a distributive lattice with the class
  represented by $\Homeo(I)$ as a global minimum.  The map
  \begin{equation*}
    C:(\Rep(I)/_{\Homeo(I)},\vee,\wedge)\to (\Pot_c(I),\cup,\cap)
  \end{equation*}
  from Proposition~\ref{prop:Cbij} is then an isomorphism of distributive
  lattices.
  \qed
\end{proposition}

\section{Spaces of traces}
\label{se:traces}

\subsection{The compact-open topology on path spaces}

We now turn to spaces of equivalence classes of paths in a Hausdorff space
$X$. To start with, we give a handy characterization of the compact-open
topology of the space of all paths $P(X)=X^I$ in $X$. By definition, is
has as a basis the sets $P(X)(\mathcal{C},\mathcal{V}):=\{p\in P(X)\mid
p(C_i)\subset V_i\}$, indexed by collections $\mathcal{C}=\{C_1,\dots
,C_n\}$ of compact subsets $C_i \subseteq I$,
resp.~$\mathcal{V}=\{V_1,\dots ,V_n\}$ of open subsets $V_i\subseteq X,\;
n\in \mathbf{N}$.

A partition $A=\{0=a_0<a_1<\dots<a_{n-1}<a_n=1\}$ of the unit interval
$I$ gives rise to the collection $\mathcal{A}=\{[a_{j-1},a_j]\}$ of
closed intervals. We write
$P(X)(A,\mathcal{U})=P(X)(\mathcal{A},\mathcal{U})\subseteq P(X)$.

\begin{lemma}\label{lem:CO}
 The sets $P(X)(A,\mathcal{U})\subset P(X)$ form a basis for the
 compact-open topology of $P(X)$.
\end{lemma}

\begin{proof}
  Compare \cite[Ch.~XII]{Dju:66} for the first part of this proof.
  Every path $q$ in a basis set $P(X)(\mathcal{C},\mathcal{V})$
  satisfies $C_i\subset q^{-1}(V_i)$; hence $C_i$ is covered by the
  connected components of $q^{-1}(V_i)$, a set of open disjoint
  intervals.  Finitely many of those, say the intervals $J_{ij}$,
  cover $C_i$; for each of those let $I_{ij}=[\min (J_{ij}\cap C_i),
  \max (J_{ij}\cap C_i)]$. Then $I_{ij}\subset J_{ij}\subseteq
  q^{-1}(V_i)$ whence $q(I_{ij})\subseteq V_i$; moreover,
  $C_i\subseteq \bigcup _jI_{ij}$.

  Let $A'= \bigcup _{i,j}\partial I_{ij}\subset I$ denote the finite
  subset of interval boundary points. Then, for two successive points
  $a,b\in A'$ we have: $q([a,b])\subseteq \bigcap_{ [a,b]\subseteq
    I_{ij}} V_i$; this intersection is to be interpreted as the set
  $X$ if the index set is empty. We slightly alter the partition $A'$
  and arrive at a new partition $A$ if one or several boundary points
  are both the upper and the lower boundary of some of the intervals
  $I_{ij}$: If $a$ is such an upper and a lower boundary of some of
  the intervals $I_{ij}$ and hence $q(a)\in \bigcap_{a\in I_{ij}}V_i$,
  we replace $a$ with a pair $a_-,a_+$ satisfying $a'<a_-<a<a_+<a''$
  for all $a'<a<a''$ in $A'$ and such that $q([a_-,a_+])\subseteq
  \bigcap_{a\in I_{ij}}V_i$.

  Let $0=a_0<a_1<\dots<a_{n-1}<a_n=1$ denote the elements of $A$ in
  the induced order, and let $\mathcal{U}$ denote the collection of
  open sets $U_j:=\bigcap _{q([a_{j-1},a_j])\subseteq V_i}V_i$. The
  open set $P(X)(\mathcal{C},\mathcal{V})$ from above then contains
  the open neighbourhood $P(X)(A,\mathcal{U})$ of $q$.
\end{proof}

\subsection{Regular traces versus traces}

In this section we compare several spaces of paths in a Hausdorff space
$X$ \emph{up to reparametrization}. Extending
Definition~\ref{def:stopreg}, we get
\begin{definition}
  \begin{enumerate}
  \item A path $p:I\to X$ in a topological space $X$ is said to be
  \emph{regular} if $\Delta_p=\emptyset$ or if $\Delta_p=\{I\}$.
\item The set of regular paths in $X$ is denoted $R(X)$ and regarded
  as a subspace of $P(X)=X^I$.
\item The spaces of (regular) paths $p$ starting in
  $x\in X$ and ending in $y\in X$ ($p(0)=x$, $p(1)=y$)
  are denoted by $R(X)(x,y)\subset P(X)(x,y)$; they are equipped with
  the induced topologies.
  \end{enumerate}
\end{definition}

Composition on the right yields a group action of the topological
group $\Homeo(I)$ on $R(X)$ and a monoid action of the topological
monoid $\Rep(I)$ on $P(X)$. These actions respect the decompositions in
subspaces $R(X)(x,y)$, resp.~$P(X)(x,y)$.

In Definition \ref{def:repeq}, we called paths $p,q\in \Rep(I)$
reparametrization equivalent if there exist reparametrizations $\phi,
\psi\in \Rep(I)$ such that $p\circ \phi =q\circ \psi$.

\begin{corollary}
  \label{co:repeq-eq}
  \begin{enumerate}
  \item Reparametrization equivalence of paths is an equivalence
    relation.
  \item Two reparametrization equivalent paths are thinly
    homotopic.
  \end{enumerate}
\end{corollary}

For a definition of thin homotopy see \cite{HKK-Hom2Gpd}; essentially, a
homotopy $H: I\times I\to X$ fixing the endpoints is thin if it factors
through a \emph{tree} (the geometric realisation of an acyclic
one-dimensional simplicial set), i.e.\ if $H: I\times I\to J\to X$ for a
tree $J$. Remark that reparametrization equivalent paths have the same
image: $p(I)=q(I)\subseteq X$. This is not necessarily true for thinly
homotopic paths; e.g., the cancellation homotopy \cite[p.~48]{Spanier}
between the concatenation of a path with its inverse and the constant path
is thin.

\begin{proof}
  \begin{enumerate}
  \item Reparametrization equivalence is clearly a reflexive and
    symmetric relation.  For transitivity, let $p,q,r\in P(X)$ denote
    three paths and assume that $p\circ \phi =q\circ \psi$ and $
    q\circ \phi'= r\circ \psi'$ for reparametrizations $\phi, \phi',
    \psi, \psi'\in \Rep(I)$. By Proposition \ref{prop:surjeq}, there
    are $\eta, \eta'\in \Rep(I)$ such that $\psi\circ \eta= \phi'\circ
    \eta'$; hence $p\circ \phi\circ \eta= r\circ \psi'\circ \eta'$.
  \item It is enough to show that $p$ and $p\circ \phi$ are thinly
    homotopic for every $\phi\in\Rep(I)$; consider the homotopy
    $H:I\times I\to I\stackrel{p}{\to} X,\;
    H(s,t)=p((1-s)t+s\phi(t))$, that even factors over $I$.
  \end{enumerate}
\end{proof}

Factoring out the respective equivalence relations given by the
actions above, we arrive at quotient spaces
$T_R(X)=R(X)/_{\Homeo(I)}$, resp.~$T(X)=P(X)/_{\Rep(I)}$ with
subspaces $T_R(X)(x,y)=R(X)(x,y)/_{\Homeo(I)}$, resp.
$T(X)(x,y)=P(X)(x,y)/_{\Rep(I)}$ for $x,y\in X$. They are considered as
spaces of (regular) \emph{traces} of paths in $X$ and should be
compared to the notions of curves or regular curves in elementary
differential geometry.

These spaces can be organised in a topological category $T(X)$ (and
likewise $T_R(X)$) with the elements of $X$ as objects, with the
topological spaces $T(X)(x,y)$ as morphism from $x$ to $y$ and with a
composition $T(X)(x,y)\times T(X)(y,z)\to T(X)(x,z)$ induced by
concatenation.  Remark that one does \emph{not} obtain a category
structure on $P(X)$ since concatenation is not associative ``on the
nose''. The categories $T(X)$ and their directed relatives are used as
important tools in \cite{Raussen:06}.

\begin{lemma}
  \label{lem:freeact}
  Let $x, y$ be elements of a Hausdorff space $X$ and let $p\in
  R(X)(x,y)$, $\phi\in \Rep(I)$. If $p\circ \phi=p$, then $\phi =\id_I$
  or $p$ is constant. In particular, for $x\neq y$, the action of
  $\Homeo(I)$ on $R(X)(x,y)$ is \emph{free}.
\end{lemma}

\begin{proof}
  If $\phi \neq \id_I$, there exists an interval $J=[a,b]\subseteq
  I$ with $\phi (a)=a, \phi (b)=b$ and, without loss of
  generality, $\phi (t)<t$ for all $a<t<b$. For all these $t$ we
  conclude that $p(t)=p(\phi (t))=p(\phi^n(t))$ for all $n>0$,
  and hence that $p(t)=p(\lim_{n\to \infty}(\phi^n(t)))=p(a)$. In
  particular, there is a non-trivial interval on which $p$ is
  constant; this is not allowed for a regular path unless $p$ is
  constant on the entire unit interval $I$ (and thus $x=y$).
\end{proof}

\begin{corollary}
  Let $x\neq y$ be elements of a topological space $X$. The quotient map
  $R(X)(x,y) \to T_R(X)(x,y)$ is a weak homotopy equivalence.
\end{corollary}

\begin{proof}
   The free group action yields a fibration with contractible fiber
   $\Homeo(I)$.
 \end{proof}

It is not clear to the authors whether one can sort out conditions
under which the quotient map is a genuine homotopy equivalence.

\subsection{Spaces of  (regular)  traces}
For $x,y\in X$, the inclusion map $R(X)(x,y)\hookrightarrow P(X)(x,y)$
induces a natural map $i:T_R(X)(x,y)\to T(X)(x,y)$ between the
corresponding quotient trace spaces. The main aim of this
section is a proof of

\begin{theorem}\label{prop:tracehomeo}
   For every two elements $x,y\in X$ of a Hausdorff space $X$, the map
   $i:T_R(X)(x,y)\to T(X)(x,y)$ is a \emph{homeomorphism}.
 \end{theorem}

 In particular, every trace can be represented by a regular trace
 (cf.~Proposition \ref{prop:sur} below). It turns out that many of the
 results on reparametrizations from the preceding section will be used in
 the proof. In a first step, we show that the map $i$ is surjective:

 \begin{proposition}
   \label{prop:sur}
   For every path $p\in P(X)$, there exists a regular path $q\in R(X)$ and
   a reparametrization $\phi\in \Rep(I)$ such that $p= q\circ \phi$.
 \end{proposition}

 \begin{proof}
   For every interval $J\subseteq I$ let $m(J)$ denote its midpoint.
   Let $m:\Delta _p\to I$ denote the map $J\mapsto m(J)$ with image
   $C:=m(\Delta _p)\subseteq I$. In order to arrive at a
   reparametrization with stop map $m$, we check the conditions from
   Proposition~\ref{prop:repar-const} for the order-preserving
   bijection $m:\Delta _p\to C$: Let $x_n= \max(J_n)\uparrow x\in I$
   denote a strictly monotone sequence. Then the midpoints converge
   as well: $m(J_n)\uparrow x$. Likewise for a decreasing sequence of
   lower boundaries and corresponding midpoints.  From Proposition
   \ref{prop:repar-const} we conclude that there exists a
   reparametrization $\phi\in \Rep(I)$ with $\Delta_{\phi}=\Delta_p$.
   Hence, there is a set-theoretic factorization
   \begin{equation*}
     \xymatrix{%
       I \ar[r]^p \ar[d]_{\phi} & X \\
       I \ar@{-->}[ur]_{q}
     }\end{equation*}
   through a \emph{regular} map $q:I\to X$.

   To check that $q$ is continuous, choose an open set $U\subset X$ and
   note that $q^{-1}(U)=\phi p^{-1}(U)$.  Combining Lemma \ref{lem:invar}
   and Lemma \ref{lem:pinterv}.3, we conclude that $p^{-1}(U)$ is a union
   of maximal disjoint open intervals that are all mapped onto open
   intervals under $\phi$.  Hence $q^{-1}(U)$ is open as well.
 \end{proof}

 \begin{proposition}\label{prop:inj}
   Two regular paths $p,q\in R(X)(x,y)$ that are reparametrization
   equivalent are in fact strictly reparametrization equivalent, i.e.,
   there exists a homeomorphism $\eta\in\Homeo(I)$ such that
   $q=p\circ\eta$.
 \end{proposition}

The proof of Proposition \ref{prop:inj} (and also that of Theorem
\ref{prop:tracehomeo}) makes use of

 \begin{lemma}
   \label{lem:inj}
   Let $p\in R(X)(x,y)$, $p'\in P(X)(x,y)$, $\phi, \phi '\in \Rep(I)$ with
   $p\circ \phi=p'\circ \phi'$. Then there exists $\eta \in \Rep(I)$ with
   $p\circ \eta =p'$. Unless $p$ is constant, $\eta$ is unique.
 \end{lemma}

 \begin{proof}
   We apply Proposition~\ref{prop:lift-left} to prove the \emph{existence}
   of such a reparametrization $\eta$: For every interval $J'\in \Delta
   _{\phi'}$, there is a unique $J\in \Delta _{p'\circ \phi'} = \Delta
   _{p\circ \phi}=\Delta _{\phi}$ such that $J'\subseteq J$.  Hence, there
   exists a unique $\eta \in \Rep(I)$ with $\phi = \eta\circ\phi'$, whence
   $p'\circ \phi'=p\circ \phi = (p\circ \eta)\circ \phi'$. Since $\phi'$
   is onto, we conclude that $p'=p\circ \eta$.

   To prove \emph{uniqueness} of the factor $\eta$, suppose that $\eta _1,
   \eta _2\in \Rep(I)$ with $p\circ \eta _1=p\circ\eta _2$. If $\eta
   _1\neq \eta_2$, one can choose an interval $J=[a,b]$ such that $\eta
   _1(a)=\eta _2(a), \eta _1(b)=\eta _2(b)$ and $\eta _1(t)<\eta _2(t)$
   for $a<t<b$ (or vice versa). Given $a<t_0=t_0'<b$, choose an increasing
   sequence $t_i$ and a decreasing sequence $t_i'$ such that
   $\eta_1(t_{i+1})=\eta_2(t_i)$, resp.~$\eta_2(t_{i+1}')=\eta_1(t_i)$.
   Both sequences converge to, say, $a\le T'<T\le b$. Since
   $\eta_1(T)=\eta_2(T)$ and $\eta_1(T')=\eta_2(T')$, we must have $T'=a,
   T=b$.

   On the other hand, $p(\eta_2(t_i))= p(\eta_1(t_{i+1}))=
   p(\eta_2(t_{i+1}))$ and hence $p(t_0)=p(b)$, and, by a similar
   argument: $p(t_0)=p(a)$. Since $t_0\in \oil a,b\oir$ was chosen
   arbitrarily, $p$ has to be constant on the interval $J=[a,b]$; this
   is impossible for a regular path $p$ unless it is constant.
 \end{proof}

 \begin{proof}(of Proposition \ref{prop:inj})
   Let $p,q\in R(X)(x,y)$ be reparametrization equivalent, i.e., there
   exist reparametrizations $\phi, \psi\in \Rep(I)$ such that $p\circ
   \phi=q\circ \psi$.  If one of $p,q$ is constant, the other is as
   well.  Assume that neither of them is constant.  By Lemma
   \ref{lem:inj}, there is a reparametrization $\rho\in \Rep(I)$ with
   $q=p\circ\rho$. Since $q$ is regular and non-constant, $\rho$ has
   to be injective and thus a homeomorphism; in particular, $p$ and
   $q$ represent the same element in $T_R(X)(x,y)$.
 \end{proof}

 \begin{proof}(of Theorem~\ref{prop:tracehomeo})
   From Propositions~\ref{prop:sur} and \ref{prop:inj} we conclude that
   the map $i$ from Theorem~\ref{prop:tracehomeo} is a continuous
   bijection. To see that it is also open, the following diagram will be
   useful:
   \begin{equation}\label{eq:tracehomeo}
     \vcenter{\xymatrix{%
         R(X)(x,y)\ar[r]^{\subseteq}\ar[d]_{Q_R} & P(X)(x,y)\ar[d]^Q\\
         T_R(X)(x,y)\ar[r]^i & T(X)(x,y)
       }}
   \end{equation}

   By Lemma~\ref{lem:CO}, a basis for the topology of $R(X)( x,y)$ is
   given by the sets
  \begin{equation*}
    R(X)( A, \mathcal U; x,y)= \{ p\in R(X)( x,y)\mid p([ a_{j-1},
    a_j])\subseteq U_j\}
  \end{equation*}
  indexed by all partitions $A=\{ 0= a_0<\dots< a_n= 1\}$ of $I$ and
  collections of $n$ open sets $\mathcal U=\{ U_1,\dots, U_n\},\;
  n\in\mathbf{N}$.

  The sets $Q_R( R(X)( A, \mathcal U; x,y))$ form a basis for the
  quotient topology on $T_R(X)( x,y)$ since
  \begin{align*}
    Q_R^{-1}( Q_R( R(X)( A, \mathcal{U}; x,y))) &=
    \bigcup\nolimits_{\rho\in\Homeo
      (I)}\rho\cdot R(X)( A, \mathcal{U}; x,y))\\
    &= \bigcup\nolimits_{\rho\in\Homeo (I)}R(X)( \rho^{-1}(A),
    \mathcal{U}; x,y)),
  \end{align*}
  and the latter set is open in $R(X)( x,y)$.

  We need to show that the sets $i( Q_R( R(X)( A, \mathcal U; x,y)))=
  Q( R(X)( A, \mathcal U; x,y))$ are open in $T(X)( x,y)$, i.e., that
  the sets $Q^{-1}( Q( R(X)( A, \mathcal U; x,y)))$ are open in $P(X)(
  x,y)$. To this end, we note first that $Q^{-1}( Q( R(X)( A, \mathcal
  U; x,y)))= \Rep (I)\cdot R(X)( A, \mathcal{U}; x,y)$: The inclusion
  $\supseteq$ is obvious. If $p\in R(X)( A, \mathcal U; x,y)$ and
  $p'\in P(X)( x,y)$, $\phi, \phi'\in \Rep(I)$ are such that $p\circ
  \phi= p'\circ \phi'$, then Lemma~\ref{lem:inj} yields $\eta\in
  \Rep(I)$ such that $p'= p\circ \eta$, hence $p'\in \Rep (I)\cdot
  R(X)( A, \mathcal{U}; x,y)$.

  In the next two steps we show that every element $q\in \Rep (I)\cdot
  R(X)( A, \mathcal{U}; x,y)$ has an open neighbourhood in that set.
  For elements in $R(X)( A, \mathcal{U}; x,y)$ it suffices to check
  that $P(X)( A, \mathcal{U}; x,y)\subseteq \Rep(I)\cdot R(X)( A,
  \mathcal{U}; x,y)$: According to Proposition~\ref{prop:sur}, a path
  $p\in P(X)(A,\mathcal{U};x,y)$ can be factored in the form $p=q\circ
  \phi$ with $q\in R(X)(x,y), \phi\in\Rep(I)$; the factors then have to
  satisfy $\phi ([a_{i-1},a_i])\subseteq q^{-1}(U_i)$.  By
  Lemma~\ref{lem:dense}, there is an approximation $\rho\in\Homeo (I)$
  to $\phi$ satisfying $\rho ([a_{i-1},a_i])\subseteq q^{-1}(U_i)$ for
  all $i$. We can thus factor $p=(q\circ \rho)\circ (\rho^{-1}\circ
  \phi)$ with $q\circ \rho \in R(X)(A,\mathcal{U};x,y)$, and hence $p\in
  \Rep(I)\cdot R(X)(A,\mathcal{U};x,y)$.

  Now let $q\in R(X)( A, \mathcal{U}; x,y)$ and $\phi\in \Rep(I)$. Let
  $B=\{ b_0= 0<\dots< b_n=1\}$ denote any partition of $I$ with $\phi
  (B)=A$. Then $q\circ \phi\in P(X)( B, \mathcal{U}; x,y)$, which is
  open in $P(X)( x,y)$. Hence it will suffice to show that $P(X)( B,
  \mathcal{U}; x,y)\subseteq \Rep (I)\cdot R(X)( A, \mathcal{U};
  x,y)$. To this end, choose $\rho\in \Homeo(I)$ such that $\rho(
  b_i)= a_i$. Let $p\in P(X)( B, \mathcal{U}; x,y)$ and write $p=
  p_1\circ \rho$. Then $p_1\in P(X)(A,\mathcal{U};x,y)\subseteq
  \Rep(I)\cdot R(X)( A, \mathcal{U}; x,y)$ by the previous step, and
  hence also $p\in \Rep(I)\cdot R(X)( A, \mathcal{U}; x,y)$.
\end{proof}

\begin{remark}\label{rem:homtype}
  All four spaces in (\ref{eq:tracehomeo}) are (at least) weakly
  homotopy equivalent: If $x,y$ are not in the same path component,
  they are all empty. Otherwise, $P(X)(x,y)$ is homotopy equivalent to
  the loop space $\Omega(X)(x)$ based at $x$; likewise, $R(X)(x,y)$ is
  homotopy equivalent to the space of regular loops $\Omega_R(X)(x)$
  based at $x$.  Both loop spaces are fibres in fibrations
  $P_R(X)\downarrow X$, resp.~$P(X)\downarrow X$ over $X$ with
  contractible total spaces $P_R(X;x)$, resp.~$P(X;x)$ of (regular)
  paths starting at $x$. The Five Lemma shows that the inclusion map
  $\Omega_R(X)(x)\hookrightarrow \Omega(X)(x)$ is a weak homotopy
  equivalence.
\end{remark}

Note that Lemma~\ref{lem:inj} allows us to give the following
``backwards'' characterization of reparametrization equivalence:

\begin{proposition}
  \label{pr:repeq-back}
  Paths $p,q\in P(X)$ are reparametrization equivalent if and only if
  there exists $r\in P(X)$ and $\phi, \psi\in \Rep(I)$ such that $p=
  r\circ \phi$ and $q= r\circ \psi$.
\end{proposition}

\begin{proof}
  For the ``only if'' part, suppose $p\circ \phi_1=q\circ \psi_1,
  \phi_1, \psi_1\in \Rep(I)$. We use
  Proposition~\ref{prop:sur} to write $p= r_1\circ \eta_1$ and $q=
  r_2\circ \eta_2$, with $r_1, r_2\in R(X)$ and $\eta_1,
  \eta_2\in \Rep(I)$. Then $r_1\circ \eta_1\circ \phi_1= r_2\circ
  \eta_2\circ \psi_1$; by Lemma \ref{lem:inj}, there exists $\eta\in
  \Rep(I)$ such that $r_2= r_1\circ \eta$, whence $q= r_1\circ
  \eta\circ \eta_2$.

  The reverse implication is clear by Proposition~\ref{prop:surjeq};
  if $\phi_2, \psi_2\in \Rep(I)$ are such that $\phi\circ \phi_2=
  \psi\circ \psi_2$, then $p\circ \phi_2=
  r\circ\phi\circ\phi_2=r\circ\psi\circ \psi_2=q\circ \psi_2$.
\end{proof}

Let us finally note the following consequence of Proposition
\ref{prop:sur}:

\begin{definition}\label{def:lf}
  A path $p:I\to X$ is called \emph{loop-free} if $p(s)= p(t)$ for any $s<
  t\in I$ implies that the restriction $p\rest{[s,t]}$ is the constant
  path.
\end{definition}

Note that a loop-free \emph{regular} path is either constant or injective.

 \begin{corollary}\label{cor:im-I}
   A loop-free path $p: I\to X$ in a Hausdorff space $X$ has an image
   $p(I)\subseteq X$ that is either a point or homeomorphic to
   $I$.
 \end{corollary}

 \begin{proof}
   By Proposition \ref{prop:sur}, there is a factorization $p= q\circ
   \phi$ with $q$ a loop-free regular and thus either constant or
   injective path. In the second case, $q$ is a continuous bijection
   from the compact space $I$ to its image $q(I)= p(I)\subseteq X$.
   The claim follows since $X$ is Hausdorff.
 \end{proof}

\section{Directed traces}\label{s:dirtra}

Originally motivated by models arising in concurrency theory in
theoretical computer science, an investigation of topological spaces
with ``preferred directions'' has been launched under the title
``Directed Algebraic Topology''. Various frameworks (local
po-spaces~\cite{FGR-ATC}, d-spaces~\cite{Grandis-DHTP},
flows~\cite{Gaucher:03} and others) have been suggested, all modifying
in various ways concepts from elementary algebraic topology, in
particular replacing relevant groups or groupoids by categories. The
d-spaces introduced by M.~Grandis \cite{Grandis-DHTP} have turned out
to give rise to a particularly successful way to combine homotopy
theoretical and categorical methods for the study of spaces with
preferred directions.

The question whether one can neglect reparametrizations in a homotopy
theoretical study of spaces of directed paths was one of the original
motivations for this article. This is -- under a certain natural
condition -- affirmed in Corollary \ref{cor:dsur}, which is used as a
starting point for the homotopy theoretical study of ``directed spaces''
in \cite{Raussen:06}.

\begin{definition}[\cite{Grandis-DHTP}]\label{def:d}
  A \emph{d-space} is a topological space $X$ together with a set
  $\vec{P}(X)\subseteq P(X)$ of continuous paths $I\to X$ such that
  \begin{enumerate}
  \item $\vec{P}(X)$ contains all constant paths;
  \item $p\circ \phi\in \vec{P}(X)$ for any $p\in \vec{P}(X)$ and any
    continuous increasing (not necessarily surjective)s map $\phi: I\to
    I$;
  \item for all $p,q\in \vec{P}(X)$ such that $p(1)= q(0)$, their
    concatenation $p*q\in \vec{P}(X)$, cf.~(\ref{eq:concat}).
  \end{enumerate}
  Elements of $\vec{P}(X)$ are called d-\emph{paths}.  $\vec{P}(X)\subseteq
  P(X)$ is given the subspace topology (of the compact-open
  topology).
\end{definition}

\begin{definition}\label{def:dmap}
  A \emph{d-map} between d-spaces $X$, $Y$ is a continuous mapping $f:
  X\to Y$ satisfying $p\in \vec{P}(X)\Rightarrow f\circ p\in \vec{P}(Y)$.
  Isomorphisms in the category of d-spaces are called
  \emph{d-homeomorphisms}.
\end{definition}

The d-interval $\vec{I}$ is given the standard d-structure
$\vec{P}(I)=\Rep(I)$.

In general d-spaces, it may occur that a non-d-path becomes directed after
reparametrization. To exclude this possibility, we add

\begin{definition}\label{def:sat}
  A d-space $X$ is called \emph{saturated} if it has the following
  additional property:
  \begin{enumerate}
  \item[4.] If $p\in P(X), \phi \in \Rep(I)$ and $p\circ \phi\in \vec{P}(X)$,
    then $p\in \vec{P}(X)$.
  \end{enumerate}
\end{definition}

In words: If a path becomes a d-path after reparametrization, then it has
to be a d-path itself already. Remark that for a saturated d-space, two
trace equivalent paths are either both d-paths, or neither of them is. The
d-space $\vec{I}$ is saturated. It is easy to turn a given d-space $X$
into a saturated one: one just adds all paths $p\in P(X)$ for which there
is a reparametrization $\phi \in \Rep(I)$ with $p\circ \phi\in \vec{P}(X)$ to
the d-paths in a new structure $S\vec{P}(X)$, which is easily seen to satisfy
the properties of a saturated d-space. So there is no harm in assuming
that a d-space is saturated right away.

Among the d-paths in $X$, we pay particular attention to the
\emph{regular} d-paths, cf.~Definition \ref{def:stopreg}.3; the set of
all those will be denoted $\vec{R}(X)=R(X)\cap \vec{P}(X)$ and equipped with
the subspace topology. Again, $\Homeo (I)$ acts (essentially freely)
on $\vec{R}(X)$, and $\Rep(I)$ acts on $\vec{P}(X)$.  We can now speak of
spaces of (regular) traces in a saturated d-space $X$:

\begin{definition}
  \begin{itemize}
  \item $\vec{T}_R(X)(x,y):=\vec{R}(X)(x,y)/_{\Homeo(I)}\subseteq T_R(X)(x,y)$
  \item $\vec{T}(X)(x,y):=\vec{P}(X)(x,y)/_{\Rep(I)}\subseteq T(X)(x,y)$
  \end{itemize}
\end{definition}

These form the morphisms of categories $\vec{T}_R(X)$, resp.~$\vec{T}(X)$ that
are investigated from a homotopy theory point of view in
\cite{Raussen:06}. The following consequence of Theorem
\ref{prop:tracehomeo} tells us that it makes no difference in topology
which of the two (quotient) trace spaces is chosen:

\begin{corollary} \label{cor:dsur} Let $X$ denote a saturated d-space
and let $x,y\in X$. The map $\vec{i}:\vec{T}_R(X)(x,y)\to \vec{T}(X)(x,y)$
induced by inclusion $\vec{R}(X)(x,y)\hookrightarrow \vec{P}(X)(x,y)$ is a
homeomorphism.
\end{corollary}

\begin{remark}
  It is no longer clear whether the inclusion map
  $\vec{R}(X)(x,y)\hookrightarrow \vec{P}(X)(x,y)$ is a weak homotopy
  equivalence. The (weak) homotopy types of both spaces depend on the
  choice of $x$ and $y$ and it is therefore not possible to argue using
  loop spaces as in \ref{rem:homtype}. From the diagram
  \begin{equation*}
    \xymatrix{%
      \vec{R}(X)(x,y)\ar[r]^{\subseteq}\ar[d]_{Q_R}^{\simeq} &
      \vec{P}(X)(x,y)\ar[d]^Q\\
      \vec{T}_R(X)(x,y)\ar[r]^i_{\cong}  & \vec{T}(X)(x,y),}
  \end{equation*}
  obtained from (\ref{eq:tracehomeo}) by restricting to d-paths, we can
  only deduce that the inclusion maps $\vec{R}(X)(x,y)\hookrightarrow
  \vec{P}(X)(x,y)$ induce injections and that the quotient maps
  $Q:\vec{P}(X)(x,y)\to \vec{T}(X)(x,y)$ induce surjections on all
  homotopy groups.
\end{remark}

\begin{definition}\label{def:d-hom}
  \begin{enumerate}
  \item \textup{\cite{Grandis-DHTP}} A \emph{d-homotopy} from a d-path
    $p\in \vec{P}(X)$ to a d-path $q\in \vec{P}(X)$ is a d-map $H:
    \vec{I}\times \vec{I}\to X$ for which $H(0,\cdot)= p$, $H(1,\cdot)=
    q$, and $H(\cdot,0)$, $H(\cdot,1)$ are constant.
  \item A d-homotopy is said to be \emph{thin} if it factors through
    the d-interval $\vec{I}$, i.e.\ if there are d-maps $\Phi: \vec{I}\times
    \vec{I}\to \vec{I}$, $r: \vec{I}\to X$ such that $H= r\circ \Phi$.
  \item Two d-paths $p,q\in \vec{P}(X)$ are said to be \emph{d-homotopic},
    respectively \emph{thinly d-homotopic}, if there exists a sequence
    $H_1,\dots, H_{2n+1}$ of d-homotopies, respectively thin d-homotopies,
    such that $H_1(0,\cdot)= p$, $H_{2n+1}(1,\cdot)= q$,
    $H_{2i-1}(1,\cdot)= H_{2i}(1,\cdot)$, and $H_{2i}(0,\cdot)=
    H_{2i+1}(0,\cdot)$.
  \end{enumerate}
\end{definition}

\begin{remark}
  \begin{itemize}
  \item The directed structure on a product of d-spaces $X$ and $Y$ is
    given by
    \begin{equation*}
      P(X\times Y)\supseteq \vec{P}(X\times Y)=\vec{P}(X)\times
      \vec{P}(Y)\subseteq P(X)\times P(Y)
    \end{equation*}
    under the natural identification $P(X\times Y)\cong P(X)\times P(Y)$.
    In particular, $\vec{P}( \vec{I}\times \vec{I})$ consists of all paths $p: I\to
    I\times I$ that are (weakly) increasing in both coordinates.
  \item The relations $\preccurlyeq$, $\preccurlyeq_T$ on d-paths
    given by existence of (thin) d-homotopies are preorders on
    $\vec{P}(X)$. The relations $\simeq$, $\simeq_T$ on d-paths given by
    being (thinly) d-homotopic are equivalence relations on $\vec{P}(X)$;
    they are the symmetric, transitive closures of $\preccurlyeq$
    respectively $\preccurlyeq_T$.
  \item In the differentiable setting, thin homotopies of
    differentiable loops and paths were defined (using homotopies of
    rank at most 1) and a resulting smooth fundamental groupoid used
    for holonomy considerations in \cite{CP:94,MP:02}.
  \end{itemize}
\end{remark}

\begin{proposition}
  \label{pr:traeqdhom}
  Two d-paths $p, q$ in a saturated d-space are reparametrization
  equivalent if and only if they are thinly d-homotopic.
\end{proposition}

This result is needed in the study \cite{Fahrenberg:06} of directed
squares. Note that the notion of d-homotopy factors over $\vec{T}(X)$ (and
$\vec{T}_R(X)$).

\begin{proof}
  We use Proposition~\ref{pr:repeq-back}. For the forward implication,
  write $p= r\circ \phi$, $q= r\circ \psi$, and let $\eta= \max(
  \phi, \psi)$. Define $\Phi, \Psi: \vec{I}\times \vec{I}\to \vec{I}$ by
  \begin{align*}
    \Phi(s,t) &= (1-s) \phi(t)+ s \eta(t)\\
    \Psi(s,t) &= (1-s) \psi(t)+ s \eta(t)
  \end{align*}
  then $r\circ \Phi$, $r\circ \Psi$ are thin d-homotopies connecting
  $p$ and $q$.

  To show the back implication, it is enough to consider the case
  where $p$ and $q$ are connected by \emph{one} thin d-homotopy $H$
  with $H(0,\cdot)= p$ and $H(1,\cdot)= q$. Write $H= r\circ \Phi:
  \vec{I}\times \vec{I}\to \vec{I}\to X$; by Corollary~\ref{cor:dsur},
  we can assume $r$ to be regular. Also, by reparametrizing if
  necessary, we can assume that $\Phi(0,0)= 0$ and $\Phi(1,1)= 1$.
  Then
  \begin{align*}
    r( \Phi(s,0)) &= H(s,0)= H(0,0)= r(0) \\
    r( \Phi(s,1)) &= H(s,1)= H(1,1)= r(1)
  \end{align*}
  hence by regularity of $r$ and following the logic of the final step in
  the proof of Lemma~\ref{lem:freeact}, $\Phi(s,0)= 0$ and $\Phi(s,1)= 1$.

  Now define $\phi, \psi: \vec{I}\to \vec{I}$ by $\phi(t)= \Phi(0,t)$,
  $\psi(t)= \Phi(1,t)$, then $\phi, \psi\in \Rep(I)$ and $p= r\circ
  \phi$, $q= r\circ \psi$.
\end{proof}

Finally, we modify the results from Section \ref{se:traces} about
loop-free paths (Definition \ref{def:lf}) to the d-space environment:
\begin{definition}
  A d-space $X$ is said to be \emph{locally loop-free} provided that
  every point has a neighbourhood in which all non-constant d-paths
  are loop-free.
\end{definition}

Note that d-spaces arising from a space with a locally partial
order~\cite{FGR-ATC} are locally loop-free. The following result
applies such to such spaces, in particular. For po-spaces, a result
similar to the following had previously been obtained in
\cite[Thm.~5]{Haucourt:05}.

\begin{corollary}\label{cor:diint}
  If $p\in \vec{P}(X)$ is a loop-free d-path in a locally loop-free
  saturated Hausdorff d-space $X$, then its image $p(\vec{I})$ is either a
  point or d-homeomorphic to $\vec{I}$.
\end{corollary}

\begin{proof}
  The statement is trivial for a constant d-path. Otherwise,
  Corollary~\ref{cor:dsur} provides us with a \emph{regular} loop-free
  d-path $q: \vec{I}\to X$ with $p=q\circ \phi$, $\phi \in \Rep(I)$ which
  by Corollary~\ref{cor:im-I} yields the homeomorphism $q:I\to
  p(\vec{I})$.  All we need to show is that its inverse $q^{-1}:
  p(\vec{I})\to \vec{I}$ is a d-map.

  Let $r: \vec{I}\to p(\vec{I})$ be a d-path; we need to show that
  $q^{-1}\circ r\in \vec{P}(I)= \Rep(I)$. Let $t_1< t_2\in I$ and suppose
  that $q^{-1}( r( t_1))> q^{-1}( r( t_2))$.

  Restricting to a smaller interval, if necessary, will ensure that
  $r([t_1, t_2])\subseteq X$ is contained in a loop-free neighbourhood
  $U\subseteq X$. The concatenation
  \begin{equation*}
    r\rest{[t_1, t_2]}* q\rest{[q^{-1}( r( t_2)), q^{-1}( r( t_1))]}
  \end{equation*}
  is a d-path and a loop in $U$ and hence constant. Then
  $q\rest{[q^{-1}( r( t_2)), q^{-1}( r( t_1))]}$ is constant, in
  contradiction to being a regular and non-constant path. Hence
  $q^{-1}( r( t_1))\le q^{-1}( r( t_2))$ whence $q^{-1}$ is a d-map.
\end{proof}

\begin{remark}
  It seems plausible that many of the methods and of the results from this
  article allow generalizations to maps $p:I^n\to X$, resp.\
  $p:\vec{I}^n\to X$, from (directed) cubes to (d)-spaces. The relevant
  reparametrizations to investigate are the d-maps $\phi:\vec{I}^n\to
  \vec{I}^n$ (monotone in every coordinate) that preserve boundaries in
  the following sense:
  \begin{equation*}
    \phi (x_1,\dots,x_n)=(y_1,\dots,y_n) \mbox{ and } x_i=0, \mbox{
      resp. }1\Rightarrow y_i=0, \mbox{ resp. }1.
  \end{equation*}
  In the differentiable setting, such reparametrizations for so-called
  $n$-loops and resulting smooth homotopy groups have been defined and
  studied in \cite{CP:98,MP:02}.
\end{remark}

\newcommand{\Afirst}[1]{#1}

\end{document}